\documentclass{article}
\usepackage[latin1]{inputenc}
\usepackage[T1]{fontenc} 

\usepackage{graphicx}
\usepackage{subfigure}          

\usepackage[english]{babel}
\usepackage[autostyle=true,german=quotes]{csquotes}
\usepackage{ngerman} 

\usepackage{placeins}

\usepackage{tabularx}

\usepackage{amsmath}
\usepackage{amsfonts}
\usepackage{amssymb}

\usepackage{lmodern,textcomp}

\usepackage{algorithm}
\usepackage{algpseudocode}

\usepackage{enumitem}
\newlist{Aufz}{enumerate}{10}
\setlist[Aufz]{label*=\arabic*.}

\usepackage{xspace}
\newcommand\largeparbreak{\par\bigskip}

\usepackage{ntheorem}

\newcommand{\fuer}{\,,\qquad}

\newcommand{\R}{\mathbb{R}\xspace}
\newcommand{\N}{\mathbb{N}\xspace}
\newcommand{\C}{\mathbb{C}\xspace}

\newcommand{\cV}{\mathcal{V}\xspace}

\newcommand{\cP}{\mathcal{P}\xspace}

\newcommand{\cG}{\mathcal{G}\xspace}

\newcommand{\cJ}{\mathcal{J}\xspace}
\newcommand{\cU}{\mathcal{U}\xspace}

\newcommand{\cK}{\mathcal{K}\xspace}
\newcommand{\cS}{\mathcal{S}\xspace}

\newcommand{\cO}{\mathcal{O}\xspace}

\newcommand{\CG}{\textsf{CG}\xspace}

\newcommand{\IDR}{\textsf{IDR($s$)}\xspace}
\newcommand{\IDReins}{\textsf{IDR($1$)}\xspace}
\newcommand{\SRGMRES}{\textsf{SRGMRES}\xspace}
\newcommand{\IDRO}{\textsf{IDR}\xspace}

\newcommand{\SRIDR}{\textsf{SRIDR($s$)}\xspace}
\newcommand{\SRIDRO}{\textsf{SRIDR}\xspace}
\newcommand{\SRBiCG}{\textsf{SRBiCG}\xspace}
\newcommand{\SRCG}{\textsf{SRCG}\xspace}

\newcommand{\BiCG}{\textsf{BiCG}\xspace}
\newcommand{\RBiCG}{\textsf{RBiCG}\xspace}
\newcommand{\RGCR}{\textsf{RGCR}\xspace}

\newcommand{\GMRES}{\textsf{GMRES}\xspace}
\newcommand{\FOM}{\textsf{FOM}\xspace}

\newcommand{\GMRESo}{\textsf{GMRES}}
\newcommand{\MLBiCG}{\textsf{ML($k$)BiCG}\xspace}

\newcommand{\MINRES}{\textsf{MINRES}\xspace}
\newcommand{\RMINRES}{\textsf{R-MINRES}\xspace}

\newcommand{\SRMRno}{\textsf{SRMINRES}}
\newcommand{\SRMR}{\SRMRno\xspace}
\newcommand{\SRMLBiCG}{\textsf{SRML($k$)BiCG}\xspace}

\newcommand{\MLbiL}{{ML-biL.-dec.}\xspace}
\newcommand{\rhs}{\textsf{rhs}\xspace}

\newcommand{\MV}{\textsf{MV}\xspace}
\newcommand{\rMV}{\textsf{\#MV}\xspace}

\newcommand{\rRD}{\textsf{\#RD}\xspace}

\newcommand{\bA}{\textbf{A}\xspace}

\newcommand{\bB}{\textbf{B}\xspace}

\newcommand{\hbV}{\hat{\textbf{V}}\xspace}
\newcommand{\tbW}{\tilde{\textbf{W}}\xspace}
\newcommand{\hbW}{\hat{\textbf{W}}\xspace}

\newcommand{\bU}{\textbf{U}\xspace}
\newcommand{\bV}{\textbf{V}\xspace}
\newcommand{\bZ}{\textbf{Z}\xspace}

\newcommand{\bP}{\textbf{P}\xspace}
\newcommand{\bPi}{\boldsymbol{\Pi}\xspace}

\newcommand{\bK}{\textbf{K}\xspace}
\newcommand{\bF}{\textbf{F}\xspace}
\newcommand{\bx}{\textbf{x}\xspace}

\newcommand{\by}{\textbf{y}\xspace}
\newcommand{\bb}{\textbf{b}\xspace}
\newcommand{\br}{\textbf{r}\xspace}
\newcommand{\bc}{\textbf{c}\xspace}
\newcommand{\bg}{\textbf{g}\xspace}
\newcommand{\bt}{\textbf{t}\xspace}

\newcommand{\hbw}{\hat{\textbf{w}}\xspace}

\newcommand{\bp}{\textbf{p}\xspace}
\newcommand{\tbU}{\tilde{\textbf{U}}\xspace}

\newcommand{\tbV}{\tilde{\textbf{V}}\xspace}

\newcommand{\bz}{\textbf{z}\xspace}
\newcommand{\be}{\textbf{e}\xspace}
\newcommand{\bw}{\textbf{w}\xspace}
\newcommand{\bM}{\textbf{M}\xspace}

\newcommand{\bC}{\textbf{C}\xspace}

\newcommand{\bW}{\textbf{W}\xspace}

\newcommand{\bH}{\textbf{H}\xspace}
\newcommand{\bI}{\textbf{I}\xspace}
\newcommand{\bR}{\textbf{R}\xspace}
\newcommand{\bL}{\textbf{L}\xspace}
\newcommand{\bT}{\textbf{T}\xspace}

\newcommand{\bO}{\textbf{0}\xspace}

\newcommand{\tbx}{\tilde{\textbf{x}}\xspace}
\newcommand{\tbb}{\tilde{\textbf{b}}\xspace}

\newcommand{\bu}{\textbf{u}\xspace}
\newcommand{\bv}{\textbf{v}\xspace}
\newcommand{\nEqns}{{n_{\text{Eqns}}}\xspace}

\newcommand{\opspan}{\textsl{span}\xspace}
\newcommand{\opImage}{\textsl{range}\xspace}

\newcommand{\rank}{\textsl{rank}\xspace}

\newcommand{\obT}{\overline{\bT}\xspace}
\newcommand{\obZ}{\overline{\bZ}\xspace}

\newcommand{\obH}{\overline{\bH}\xspace}

\newcommand{\obU}{\overline{\bU}\xspace}

\newcommand{\ubT}{\underline{\bT}\xspace}

\newcommand{\ubH}{\underline{\bH}\xspace}

\newcommand{\obV}{\overline{\textbf{V}}\xspace}
\newcommand{\obW}{\overline{\textbf{W}}\xspace}
\newcommand{\obC}{\overline{\textbf{C}}\xspace}

\newcommand{\cond}{\textsl{cond}\xspace}

\newcommand{\hia}{^{(\iota)}\xspace}
\newcommand{\hib}{^{(\iota+1)}\xspace}
\newcommand{\himu}{^{(\iota+\mu)}\xspace}

\newcommand{\tol}{\mathrm{tol}\xspace}

\newcommand{\mdtheorem}[2]{\newtheorem{#1}{#2}}

\mdtheorem{Definition}{Definition}
\mdtheorem{Problem}{Problem}
\mdtheorem{Bemerkung}{Bemerkung}
\mdtheorem{Methode}{Methode}
\mdtheorem{Satz}{Satz}
\mdtheorem{Theorem}{Theorem}
\mdtheorem{Prototyp}{Prototyp}
\mdtheorem{Lemma}{Lemma}
\mdtheorem{Beispiel}{Beispiel}
\mdtheorem{Vermutung}{Vermutung}

\title{Short-Recurrence and -Storage Recycling of large Krylov-Subspaces for Sequences of Linear Systems with changing Right-Hand-Sides}
\author{Martin Neuenhofen \footnote{Martin.Peter.Neuenhofen@rwth-aachen.de}}

\date{\today}
\begin{document}

\maketitle

\begin{abstract}
In this text I present a couple of new principles and thereon based iterative methods for numerical solution of sequences of systems of linear equations with fixed system matrix and changing right-hand-sides. The use of the new methods is to recycle \textit{all} subspace information that is obtained anyway in the solution process of a former system, to solve subsequent systems.

All these principles and methods are based on short recurrences and small storage requirements. The principles are based on the \IDRO-theorem \cite{IDR-report} and the Horner scheme for polynomials.

\end{abstract}

\tableofcontents

\section{Introduction}

\subsection{Problem Statement}

In this text numerical methods for solution of problem \ref{prob:2} are presented.
\begin{Problem}[]\label{prob:2}
	Let $\bA \in \C^{N \times N}$, $\bb : \N \times \C^N \rightarrow \C^N$, $\bx^{(0)} = \bO_{N \times 1}$. Find $\bx\hia,\br\hia \in \C^N$ for $\iota = 1,...,\nEqns$, such that
	\begin{align*}
	\bA \cdot \bx\hia + \br\hia = \bb\big(\iota,\bx^{(\iota-1)}\big),\quad \iota = 1,...,\nEqns\,,
	\end{align*}
	and $\|\br\hia\|$ are small.
\end{Problem}
\paragraph{Remark}
As the right-hand-sides are each dependent on the solution of the former right-hand-side (\rhs), all systems have to be solved one after another.

Problems of this structure result e.g. from the discretization of the partial differential equation (pde)
\begin{align}
	\begin{Bmatrix}
		\begin{split}
			\frac{\partial u(x,y,t)}{\partial t}-\Delta u(x,y,t) &= 1 \quad (x,y,t) \in \Omega_T:= (0,1)^2 \times \left(0,1\right] \\
			u(x,y,t)&=0 \quad (x,y,t) \in \partial\Omega_T \setminus \Omega_T
		\end{split}
	\end{Bmatrix} \label{eqn:Fourier}
\end{align}
with central finite differences with $\Delta x = 1/101$ and implicit Euler's method with $\Delta t = 0.1$ to
\begin{align*}
	\bA &= \bI + 0.1 \cdot 101^2 \cdot \hat{\bA} \in \R^{10000 \times 10000}\\
	\bb^{(\iota)} &= \bx^{(\iota-1)} + 0.1 \cdot \textbf{1} \in \R^{10000}\quad \iota = 1,2,...,10\,.
\end{align*}
with $\hat\bA$ the matrix to the 2-dimensional $5$-point-(-1,-1,4,-1,-1)-stencil or in Matlab-notation $\hat{\bA} = \texttt{gallery('poisson',100)}$. For \texttt{gallery} see \cite{MATLAB-TEST}.

\subsection{Literature Review}
A selection of approaches for numerical treatment of problem \ref{prob:2} in both Hermitian and non-Hermitian case by means of Krylov-subspace methods is given in \cite{GCRO-DR,R-MINRES,R-BiCG,RGCRO}.
\largeparbreak
The central idea of these and many other approaches from literature can be summarized with the help of algorithm \ref{algo:RGCR}.
\begin{algorithm} [H]
	\caption{Recycling for Generalized Conjugate Residual}
	\label{algo:RGCR}
	\begin{algorithmic}[1]
		\Procedure{RGCR}{$\bA,\bU,\bV,\bx,\br$}
		\State \textit{// pre-conditions: $\bA \cdot \bU = \bV$, $\bV^H \cdot \bV = \bI$}
		\State $\boldsymbol{\omega} := \bV^H \cdot \br$ \label{algo:RGCR:Z3}
		\State $\bx := \bx + \bU\cdot \boldsymbol{\omega}$, $\br := \br - \bV \cdot \boldsymbol{\omega}$
		\For{$i=1,2,...$}
			\State \textit{// extend orthonormal image basis (Gram-Schmidt procedure)}
			\State $\bu := \br$, $\bv := \bA \cdot \bu$
			\State $\boldsymbol{\gamma} := \bV^H \cdot \bv$
			\State $\bu := \bu - \bU\cdot \boldsymbol{\gamma}$, $\bv := \bv - \bV \cdot\boldsymbol{\gamma}$
			\State $\ell := \|\bv\|$
			\State $\bu := 1/\ell \cdot \bu$, $\bv := 1/\ell \cdot \bv$
			\State \textit{// orthogonalize residual on new image}
			\State $\omega := \langle \bv , \br \rangle$
			\State $\bx := \bx + \bu \cdot \omega $, $\br := \br - \bv \cdot \omega$
			\State $\bU = [\bU,\bu]$, $\bV = [\bV,\bv]$ \label{algo:GCR:L14}
			\State if $\|\br\|$ is small enough then break
		\EndFor
		\State \textit{// post-conditions: $\bA \cdot \bU = \bV$, $\bV^H \cdot \bV = \bI$}
		\State \Return $\bU,\bV,\bx,\br$
		\EndProcedure
	\end{algorithmic}
\end{algorithm}
\paragraph{Remark} The use of Gram-Schmidt procedure for orthogonalization is exemplary.

The above algorithm can be used to solve problem \ref{prob:2} as follows:
\begin{algorithmic}[1]
	\State $\bU := []$, $\bV:= []$, $\bx^{(0)}:= \bO_{N \times 1}$
	\For{$\iota = 1,...,\nEqns$}
		\State $\bb\hia \leftarrow \bb(\iota,\bx^{(\iota-1)})$
		\State $[\,\bU,\bV,\bx\hia,\br\hia\,] \leftarrow \texttt{RGCR}(\,\bA,\bU,\bV,\bO_{N \times 1},\bb\hia) $
	\EndFor
\end{algorithmic}
One keeps the basis matrices $\bU,\bV$ from the former solution process for the solution process thereafter.

The size of $\bU,\bV$ increses by one in every iteration for each system of the whole sequence. This typically results in high storage and high computational cost in the Gram-Schmidt process, but may reduce the necessary number of computed matrix-vector-products with $\bA$ (in the following shortened by \rMV), compared to solving each system separately, cf. chap. \ref{chap:Experiments} and \cite{RGCRO}.

If $\bA$ is Hermitian, then a part of the Gram-Schmidt orthogonalizations can be replaced by the Lanczos procedure, cf. \RMINRES \cite{R-MINRES,R-MINRES2}. If $\bA$ is non-Hermitian, then the corresponding orthogonalizations instead can be replaced by the bi-Lanczos procedure, cf. \RBiCG \cite{R-BiCG,R-BiCG2,R-BiCG3}. Of course then the method is no longer residual optimal, as the orthogonalization is replaced by a biorthogonalization.

The principle of reusing information from a former solution process for the current system is called recycling, Krylov-recycling, Krylov-subspace recycling or subspace recycling \cite{Compressing,PGCRO-DR,RGCRO}.
\largeparbreak

One wants to use recycling to reduce the number of needed matrix-vector-products with $\bA$; concurrently one wants to limit the additional storage cost and computational overhead which arises from the recycling.

For limitation of both memory and computation overhead, \textit{compressings} \cite[p. 1]{Compressing} have been developed: For this a reduction mechanism is added in algorithm \ref{algo:RGCR} below line \ref{algo:GCR:L14}. At a certain number of columns $m$ of both $\bU,\bV \in \C^{N \times m}$, the number of columns is reduced to $k < m$ by an update of the form
\begin{algorithmic}[1]
	\State $\bU := \bU \cdot \boldsymbol{\nu}$
	\State $\bV := \bV \cdot \boldsymbol{\nu}$\,,
\end{algorithmic}
where one tries to choose $\boldsymbol{\nu} \in \C^{m \times k}$ such that \rMV is minimized for all following \rhs-es, cf. \cite{GCRO-DR,GCR-OT}.

\paragraph{Remark}
If the system matrix is changed, then there exists the approach, cf. \cite{GCRO-DR,R-MINRES}, to add the following lines in \RGCR before line \ref{algo:RGCR:Z3}:
\begin{algorithmic}[1]
	\State $\bV := \bA \cdot \bU$
	\State $[\,\bV,\bR \,] \leftarrow \texttt{qr}(\bV,0)$ \textit{// reduced QR-decomposition}
	\State $\bU := \bU \cdot \bR^{-1}$
\end{algorithmic}
By this the basis $\bU$ of the search space can be recycled for the new system matrix and the columns of $\bV$ remain orthonormal w.r.t. each other.
Of course the usefullness of this approach depends much on the way how the system matrix changes.
\largeparbreak
We also noticed the approach in \cite{IDR-Milten} for \IDR. This diploma thesis deals with sequences of linear systems with changing matrices and \rhs-es.
The author proposes to hold the matrix of shadow vectors $\bP$ and the auxiliary vectors $\bV$, cf. to the notation in chap. \ref{chap:SonnRecChap}. Miltenberger describes this approach as \enquote{mimicing a higher Sonneveld-space}, cf. \cite[p.39]{IDR-Milten}. From numerical experiments this turns out to be efficient. However one can show that this approach does only lead to the reuse of a rather small search space of $s$ dimensions.

In contrast to this we will propose a method which does not only reuses $\bP$ and $\bV$ but also the relaxiations $\omega_1,...,\omega_J$ and by this truly recycles Sonneveld-spaces. This leads to the fundamental idea to all the methods which are presented here.

\subsection{Scope}
The methods presented from literature so far are unable to recycle high-dimensional search spaces. This is due to the following reasons:
\begin{Aufz}
	\item For each recycled space dimension \RGCR has to store two columns $\bu$,$\bv$ in the matrices $\bU,\bV$ ($\rightarrow$ large memory requirements).
	\item For computation of a new descent direction one has to orthogonalize it on all former descent directions ($\rightarrow$ long recurrences).
\end{Aufz}
Both limits the size of the recycling space.

In this text instead I present methods that can recycle very large Krylov-subspaces with small memory requirements and short recurrences.

To be more precise: A first system is solved with a usual Krlyov-subspace method in $n$ iterations, e.g. \MINRES, \CG, \BiCG, to orthogonalize the residual on a $n$-dimensional test space. By storing $\cO(k)$ column-vectors of size $N$ and computing $\cO(n/k)$ \MV-s with $\bA$, the proposed methods are able to find a solution to a following system with same system matrix, the residual of which is orthogonal on the $n$-dimensional test space of the former system.

\subsection{Notations}
We introduce the following two numbers and three spaces which play an important role in this text.

The spaces are: The search space $\cU$, in which the approximate solution is computed; the image space $\cV = \bA \cdot \cU$ of the search space; and the test space $\cP$, on which the residual is orthogonalized.

For Krylov-subspaces I use the following notations:
\begin{align*}
	\cK_n(\bA;\bb) &:= \operatornamewithlimits{span}_{i=1,...,n}\big\lbrace\bA^{i-1} \cdot \bb \big\rbrace\\
	\cK^\star_J(\bA;[\bv_1,...,\bv_s]) &:= \sum_{i=1,...,s} \cK_J(\bA;\bv_i)\\
	\cK_n(\bA;[\bv_1,...,\bv_s]) &:= \operatornamewithlimits{span}_{i=1,...,n}\big\lbrace\bA^{\lfloor(i-1)/s\rfloor} \cdot \bv_{\big( \,(i-1) \, \operatorname{mod} \, s)+1} \big\rbrace
\end{align*}

The important numbers are: \rMV is the number of all computed matrix-vector-products with $\bA$. \MV instead is the shortcut for \enquote{matrix-vector-product (with $\bA$)}. \rMV is a cost measure and shall be minimized. Optimally it holds: \rMV$=\dim(\cU)$. This means with each \MV the dimension of the search space increases by one.
\rRD (number of reduced dimensions) is the dimension of the test space, thus \rRD$=\dim(\cP)$.
\largeparbreak
Throughout this text the following special matrices are used: $\bH \in \C^{n \times n}$ denotes a Hessenberg-matrix, $\bU,\bV,\bW,\bC,\bZ \in \C^{N \times n}$ denote basis matrices. $n$ is the size of a Hessenberg-decomposition
\begin{align*}
	\bA\cdot\bV=\obV\cdot\obH\,.
\end{align*}
$\obH$ is a Hessenberg matrix consisting of $\bH$ with an additional row at the bottom. If $\bH$ resp. $\obH$ is banded, we instead write $\bT$ resp. $\obT$. Instead of an overbar, an underbar as in $\ubH$, $\ubT$ means an additional column right to $\bH$ resp. $\bT$.

$\obU,\obV,\obW,\obC,\obZ \in \C^{N \times (n+1)}$ are the basis matrices $\bU,\bV,\bW,\bC,\bZ$, added by one column to the right. $\bP \in \C^{N \times s}$ for $s<n$ is also a basis matrix. The column vectors of $\bU,\bV,\bW,\bC,\bZ,\bP$ are $\bu_i,\bv_i,\bw_i,\bc_i,\bz_i,\bp_i \in \C^N$, $i=1,...,n$.

In this text we use polynomials $\Omega_j(\cdot)$ of degree $j$ with non-zero leading coefficient and $\Omega_j(0)=1$.

We will often work on single column vectors of matrices. Given a matrix $\bB = [\bb_1,...,\bb_p] \in \C^{q \times p}$, $\bb_i \in \C^q$, $\bB(:,k:\ell:m) = [\bb_k,\bb_{k+\ell},\bb_{k+2 \cdot \ell},...,\bb_{k + \lfloor (m-k)/\ell \rfloor \cdot \ell}] \in \C^{q \times \lfloor (m-k)/\ell \rfloor}$ is a matrix that consists of several columns of $\bB$. The same notation holds for rows and combinations of rows and columns. This notation is equivalent to that of Matlab. For further details cf. \cite[accessing matrices]{MATLAB-TEST}.

\subsection{Structure}
The methods presented in this text approximate \RGCR.

These methods are able to reuse complete Krylov-subspace information (this means without any kind of \textit{compressing}), that is obtained anyway during the solution process of one system, to solve a following system. To achieve this they only use small storage consumption and short recurrences.

Thus the methods presented here differ from those presented in the literature review above because those need compressing to limit memory consumption and computation cost.

As the methods presented here distinguish themself by only using short recurrences and low memory, they are called \textit{Short Recycling} methods in the following.
\largeparbreak
In the second chapter a recycling principle based on \IDRO-mechanism is presented and transferred into a short recycling prototype algorithm, called \SRIDR. Then it is discussed, how this prototype can be modified under conservation of the underlying principle.

In the third chapter a principle called \textit{short representations} is presented. With this principle one can recycle information from any (generalized) Hessenberg decomposition by means of small storage and computational overhead. After the principle is described, it is applied to the symmetric and unsymmetric Lanczos process to derive several short recycling methods.

In the fourth chapter an approach is presented to achieve more robustness for methods based on short representations.

In the fifth chapter I present a strategy to increase accuracy of solutions that one obtained by the short recycling methods (henceforth called \textit{Short Recycling solutions}) from the chapters before by means of an iterative post iterations scheme.

The sixth chapter delivers a method for efficient computation of a multiple left bi-Lanczos decomposition. This can then also be used to derive a short recycling method. This chapter is rather meant as an outlook on what could be realized.

In the seventh chapter I show preconditioning approaches. For recycling of data from a generalized Hessenberg-decomposition it seems necessary to use linear non-changing preconditioners. As such preconditioners can be simply applied to the Hessenberg-decomposition itself, this chapter is more like a literature review on how to treat non-splitted preconditioners for Hermitian systems.

In the eighth chapter numerical experiments are presented for the methods \SRIDR and \SRBiCG, as for these an implementation already exists.

\section{Short Recycling with \IDRO-Principle}\label{chap:SonnRecChap}

In this section I derive a short recycling method from induced dimension reduction (\IDRO). Therefore at first we review the \IDRO-principle.

\subsection{IDR-Review}
\IDRO-methods build a sequence of residuals $\br$ and $s \in \N$ auxiliary vectors $\bV = [\bv_1,...,\bv_s]$ in Sonneveld-spaces $\cG_j$, $j=0,...,J$ of shrinking dimension. All details can be found in \cite{IDRweb}.

I give a very short introduction, also to clear the notation: By choosing $\bP \in \C^{N \times s}$, $\rank(\bP)=s$, the space $\cS = \lbrace \bv \in \C^N \vert \bP^H \cdot \bv = \bO \rbrace$ is defined.

The Sonneveld-spaces $\cG_j$, $j=0,...,J$ are then defined by the following recursion:
\begin{align*}
	\cG_0 &:= \C^N\\
	\cG_{j+1} &:= (\bI - \omega_{j+1} \cdot \bA) \cdot (\cG_j \cap \cS),\ \ \omega_{j+1} \in \C \setminus \lbrace 0 \rbrace, \quad j=0,...,J-1\,.
\end{align*}
$j$ is called \textit{level} of the Sonneveld-space.

The \IDRO-theorem now states, that in the canonical case the dimension of these spaces $\cG_j$ is reduced by $s$ with each level-increase, cf. \cite{IDR-biorth,IDR-Gutkn}.
\largeparbreak
A conventional \IDRO-algorithm uses $s$ auxiliary vectors as columns of a matrix $\bV \in \C^{N \times s}$ and the residual $\br_j$ in $\cG_j$, this means $\opImage(\bV) \subset \cG_j, \br_j \in \cG_j$, to construct a next residual $\br_{j+1}$ of higher level, this means $\br_{j+1} \in \cG_{j+1}$.

To do so, one computes
\begin{align}
	\br_{j+1} := (\bI - \omega_{j+1} \cdot \bA) \cdot \big(\bI - \bV \cdot (\bP^H \cdot \bV)^{-1} \cdot \bP^H\big) \cdot \br_j \label{eqn:res_Update_IDR}
\end{align}
and keeps the corresponding approximate solution vector $\bx$ for the new residual. A step like this is called \textit{induced dimension reduction step} \cite{IDR-report}, or short \textit{\IDRO-step}.

To repeat this computation, so that the residual can be restricted to higher and higher levels, the levels of the auxiliary vectors also have to be increased. Therefore one also computes dimension reduction steps for the auxiliary vectors. All in all one needs $s+1$ \MV-s for the residual and the $s$ auxiliary vectors to restrict both $\br$ and the columns of $\bV$ into a higher-level Sonneveld-space.

One computational loop of dimension reduction steps for the residual and all auxiliary vectors is called \textit{\IDRO-cycle}. To sum it up: To eliminate $s$ dimensions from the residual, one has to compute one \IDRO-cycle, requiring $s+1$ \MV-s with $\bA$.

It holds: \rRD$= s \cdot \lceil 1/(s+1) \,\cdot\, $\rMV $\rceil$.

\subsection{Basic Idea: Recycling Sonneveld-Spaces}
As we consider solution strategies for problem \ref{prob:2}, we now have to ask how one can use the \IDRO-mechanism to construct a short recycling method.
\largeparbreak
Consider we want to solve two systems
\begin{align*}
	\bA \cdot \bx\hia &= \bb\hia\\
	\bA \cdot \bx\hib &= \bb\hib
\end{align*}
one after another. With \IDR we can solve the first system $(\iota)$, e.g. for $s=8$ and $J^\star=10$ \IDRO-cycles. By this we obtain $\br,\bv_1,...,\bv_s \in \cG_{J^\star}$. For later reference we define $n = s \cdot J^\star$.

Now for solution of the next system with $\bb\hib$ we could still use \IDR, but this method would then compute \IDRO-steps for (new) auxiliary vectors to restrict them in higher Sonneveld-spaces. These computations are not necessary because we already have auxiliary vectors $\bv_1,...,\bv_s \in \cG_{J^\star}$ in a high Sonneveld-space.

So as long as $\br$ does not have level $ > J^\star$, one does not need to compute \IDRO-steps for the auxiliary vectors, by just holding the old auxiliary vectors $\bv_1,...,\bv_s \in \cG_{J^\star}$ and choosing the same data for construction of the Sonneveld-spaces (i.e. same $\omega_1,...,\omega_{J^\star}$, same $\bP$). By this one can restrict the residual into the next-level Sonneveld-space repeatedly by only one \IDR-step. Thus an \IDRO-cycle of this kind consists only of one \IDRO-step. We call such special \IDRO-cycles \textit{recycling-\IDRO-cycles}.

\begin{figure}
\centering
\includegraphics[width=1\linewidth]{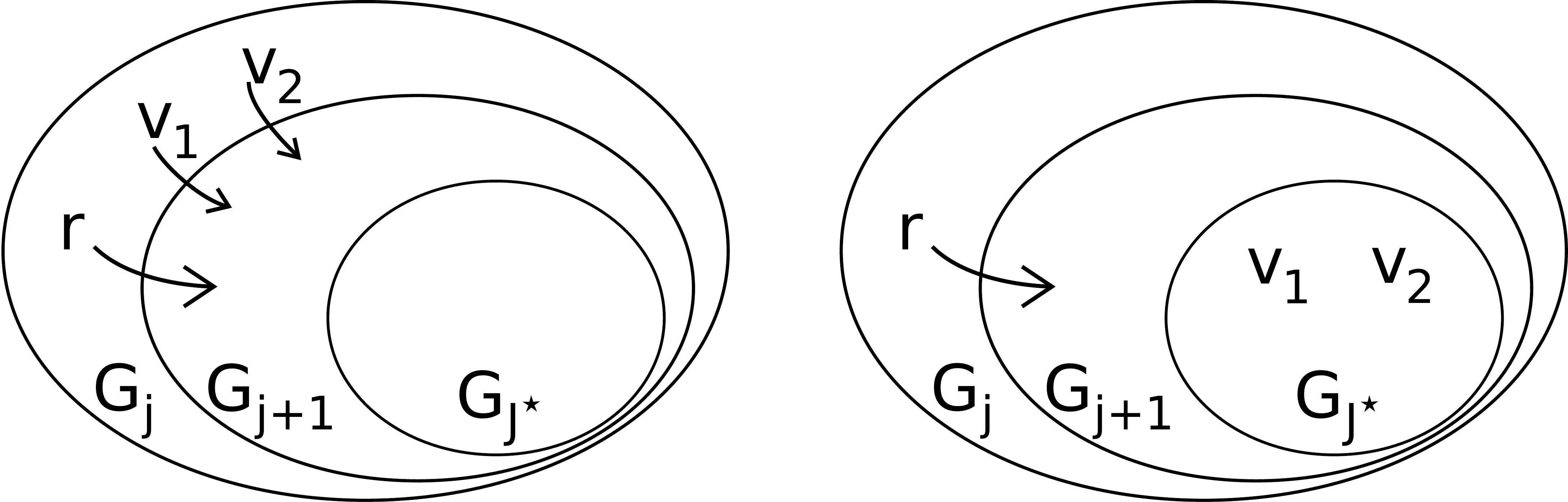}
\caption{Working principle for $s=2$. Left: Conventional \IDRO-cycles. Right: Recycling-\IDRO-cycles.}
\label{fig:SRIDR_principle}
\end{figure}

By using recycling-\IDRO-cycles, one obtains the following improved ratio between reduced dimensions and \MV-s: \rRD$= s\, \cdot\, $\rMV, as long as \rRD$\leq s \cdot J^\star$.
\largeparbreak
The basic idea of recycling-\IDRO-cycles is sketched and compared to the conventional \IDRO-cycles in fig. \ref{fig:SRIDR_principle}:

On the left side the residual and auxiliary vectors are of the same level $j$. The residual can be restricted into the next Sonneveld-space by one \IDRO-step. However to repeat this one has to perform additional \IDRO-steps to get auxiliary vectors of at least the same level as the new residual. So to repeat \IDRO-steps for the residual, one has to compute \IDRO-steps for each auxiliary vector additional to that of the residual, thus an \IDRO-cycle (i.e. a computation loop, which in repetition shrinks the residual in higher Sonneveld-space) consists of $s+1$ \IDRO-steps (i.e. $s+1$ \MV-s).

In the right part of the diagram instead, the auxiliary vectors are already in a higher-level Sonneveld-space which is predefined by the choice of $\bP$ and $\boldsymbol{\omega} = (\omega_1,...,\omega_{J^\star})$. By using these auxiliary vectors, the residual can be restricted into the next-level Sonneveld-space by an \IDRO-step without need to perform additional \IDRO-steps for the auxiliary vectors as long as $j \leq J^\star$. So as long as $j \leq J^\star$ holds one can directly afterwards compute the next \IDRO-step for the residual, thus one recycling-\IDRO-cycle consists of only $1$ \IDRO-step (i.e. $1$ \MV).

\paragraph{Original report} Best to our knowledge, the idea of \SRIDR is given the first time in \cite[part 2, chap. 1.8]{LSIDR}, published in February 2015.

\subsection{Prototype: \SRIDR}\label{chap:SRIDR_proto}
To build a prototype, that uses the idea from above, any \IDRO-algorithm can be used. The for-loop of the \IDRO-steps for the auxiliary vectors just has to be skipped if the level of $\br$ is $\leq J^\star$ and the $\omega_j$ from the solution process to the former \rhs have to be stored and reused for respective level $j \leq J^\star$. These modifiations can be found in lines \ref{algo:SRIDR_omg} and \ref{algo:SRIDR_for} of algorithm \ref{algo:SRIDR}.
\begin{algorithm} [H]
	\caption{SRIDR}
	\label{algo:SRIDR}
	\begin{algorithmic}[1]
		\Procedure{SRIDR}{$\bA,\bU,\bV,\bP,\bx,\br,\boldsymbol{\omega},J^\star,J$}
		\State \textit{// pre-conditions: $\bA \cdot \bU = \bV$, $\opImage(\bV)\subset\cG_{J^\star} $, $\rank(\bP) = s$}
		\For{$j=1,2,...,J$}
			\State $\boldsymbol{\gamma} := (\bP^H \cdot \bV)^{-1} \cdot (\bP^H \cdot \br)$
			\State $\bx := \bx + \bU \cdot \boldsymbol{\gamma}$, $\br := \br - \bV \cdot \boldsymbol{\gamma}$ \quad \textit{// $\br \in \cG_j \cap \cS$}
			\If{$j > J^\star$} \label{algo:SRIDR_omg}
			 	\State Choose $\omega_j$
			\EndIf
			\State $\bx := \bx + \omega_j \cdot \br$, $\br := \br - \omega_j \cdot \bA \cdot \br$ \quad \textit{// $\br \in \cG_{j+1}$}
			\State \textit{// If $\bA \cdot \br$ is not separately needed, e.g. because $\omega_j$ is recycled, then $\bx := \bx + \omega_j \cdot \br$, $\br := \bb-\bA \cdot \bx $ is probably better.}
			\If{$j > J^\star$} \label{algo:SRIDR_for}
				\For{$i=1,...,s$}
					\State $\boldsymbol{\gamma} := (\bP^H \cdot \bV)^{-1} \cdot (\bP^H \cdot \br)$
					\State $\bu_i := \bU \cdot \boldsymbol{\gamma} + \omega_j \cdot (\br - \bV \cdot \boldsymbol{\gamma})$
					\State $\bv_i := \bA \cdot \bu_i$ \quad \textit{// $\bv_i \in \cG_{j+1}$}
					\State \textit{// bi-orthogonalization}
					\State \textit{// ...}
				\EndFor
			\EndIf
		\EndFor
		\State \textit{// post-conditions: $\bA \cdot \bU = \bV$, $\opImage(\bV)\subset\cG_{J}$, $\rank(\bP) = s$}
		\State \Return $\bU,\bV,\bx,\br$
		\EndProcedure
	\end{algorithmic}
\end{algorithm}
We now discuss how this prototype can be improved. The modification of skipping the for-loop of \IDRO-steps for the auxiliary vectors is seen as an advantage due to reduction of \rMV, but the lost freedom in the choice of the relaxations $\omega_j$ for $j \leq J^\star$ is a drawback. Hence we now discuss remedies to this restriction.
\largeparbreak
For a subclass of problems the choice
\begin{align}
	\omega_j = {\langle\bA \cdot \br,\br\rangle}/{\|\bA \cdot \br\|^2} \label{eqn:OmegaOpt}
\end{align}
is appropriate. Therefore this is the standard choice in the algorithms in \cite{IDR-report,IDR-biorth}. So we might like to choose $\omega_j$ like this for our recycling-\IDRO-cycles but we cannot as we have to reuse a given value for $\omega_j$.
\largeparbreak
For indefinite problems the choice in \eqref{eqn:OmegaOpt} is of small use anyway. Especially for skew-Hermitian systems this choice of $\omega$ leads to a breakdown. Due to this there exist choices for $\omega$ that do not depend on the current residual but only on the system matrix. An example is
\begin{align*}
	\omega_j = 2/({\lambda_{\text{max}} + \lambda_{\text{min}}})
\end{align*}
for positive definite systems, cf. Richardson's method \cite{AMeister}.

\largeparbreak
In the following I present some other possible remedies to the choice of $\omega_j$ for $j \leq J^\star$ in algorithm \ref{algo:SRIDR}.

\paragraph{Iteration matrix} One option is to replace the \GMRESo(1)-polynomial $(\bI - \omega_j \cdot \bA)$ in the recursive definition of the $\cG_j$ by a basic linear iterative scheme like SOR or Jacobi's method \cite[chap. 4.1]{Saad1}. Then all the relaxation-parameters vanish and the Sonneveld-spaces only depend on $\bP$. A problem with this approach is that there does not necessarily exist a convergent basic iterative method for arbitrary $\bA$.

\paragraph{Change the ordering} Due to the equivalence
\begin{align}
	\cG_J = \Bigg\lbrace \bv \in \C^N \ : \ \bv = \Omega_J(\bA) \cdot \boldsymbol{\xi},\ \boldsymbol{\xi} \perp \cK^\star_{J}(\bA^H;[\bp_1,...,\bp_s]) \Bigg\rbrace \label{eqn:Sonneveld_Gut}
\end{align}
\cite[eqn. (3.11)]{IDR-Gutkn} with $\Omega_J(t) = (1-\omega_1 \cdot t) \cdot ... \cdot (1-\omega_J \cdot t)$ the order of the relaxation coefficients $\omega_j$ can be changed without destruction of the property $\bv_1,...,\bv_s \in \cG_{J^\star}$.

\paragraph{Alternating relaxations} If a relaxation $\omega$ is choosen arbitrarily for one \IDRO-cycle, then the residual of the current Sonneveld-space $\cG_j$ is introduced into a higher subspace $\tilde{\cG}_{j+1} \subset \cG_j$ because of the IDR-theorem. Thereby it remains in $\cG_j$. According to that one could for example choose $\omega$ alternatingly arbitrary or recycled. By this the ratio deteriorates from \rRD$= s \cdot $\rMV to \rRD$= s \cdot \lceil 1/2 \cdot $\rMV$\rceil$, but possibly the residual norms become smaller.
\largeparbreak

In summary the prototype offers many options for modifications. As a concluding remark, I want to give two suggestions:

\paragraph{Throw columns} From \eqref{eqn:res_Update_IDR} we see that for strictly monotonic decreasing residual norms the condition $\|I-\omega \cdot \bA\| \cdot \|\bI - \bV \cdot (\bP^H \cdot \bV)^{-1} \cdot \bP^H \| < 1$ is sufficient. $\|\bI-\omega \cdot \bA\| =:\rho < 1$ can be ensured for any positive definite system matrix $\bA$ by choosing $\omega$ small enough. $\|\bI - \bV \cdot (\bP^H \cdot \bV)^{-1} \cdot \bP^H \| < 1/\rho$ can be obtained if $\bV,\bP$ are choosen such that the canonical angles between $\opImage(\bV)$ and $\opImage(\bP)$ are sufficiently small.

To achieve this one can reduce the number of columns of $\bP,\bV \in \C^{N \times s}$ by
\begin{align*}
	\bP &:= \bP\cdot \boldsymbol{\nu}\\
	\bP &:= \bP\cdot \boldsymbol{\eta}\,,
\end{align*}
where $\boldsymbol{\nu},\boldsymbol{\eta} \in \C^{s \times \check{s}}$, $\check{s}<s$
because by reduction of $\bP$ the subspaces $\cS$ grows and therefore the old Sonneveld-spaces become subspaces of the new ones defined by the overwritten $\bP$. Thus this reduction does not destroy $\opImage(\bV) \subset \cG_{J^\star}$ but can decrease the canonical angles between $\opImage(\bV)$ and $\opImage(\bP)$.

\paragraph{\SRIDR as inner Krylov-subspace method} From \eqref{eqn:res_Update_IDR} one can see that the recycling-\IDRO-cycles of \SRIDR are linear. The reason for this is that $\bV$ does not change. Several \SRIDRO-cycles could be used for an inner linear subsolver with an outer Krylov-subspace method. This could have several uses: For example one could try to make the recycling more robust or use \SRIDR in an accelerated iterative refinement approach for sequences with slightly changing system matrices.

\subsection{Projection Spaces of \SRIDR}\label{chap:SRIDR-Spaces}
\SRIDR recycles Sonneveld-spaces. Precisely it recycles the following test space:
\begin{align*}
\cP = \Omega_{J^\star}(\bA^H) \cdot \cK^\star_{J^\star}(\bA^H;\bP)\,,
\end{align*}
where $\Omega_J^\star(t) = \prod_{j=1}^{J^\star} (1-\omega_j \cdot t)$, cf. to \eqref{eqn:Sonneveld_Gut}.

The \textit{original search space} $\cK_n(\bA;\bb\hia)$ from the first system however is not recycled.

This can be unfavorable in case one wants to make use of the knowledge that for a later system there exists a good approximate solution in the original Krylov-subspace $\cK_n(\bA;\bb\hia)$. We even see in the numerical experiments from the results of \RGCR, that this original Krylov-subspace is likely to be a very good recycling space for certain applications.

So reusing the original search space may be favorable. However for an unknown runtime and implementation dependent search space such as of \SRIDR it is much harder to give any guarantee on the usefulness of this search space. Thus we would prefer a method that recycles the original Krylov-subspace $\cK_n(\bA;\bb\hia)$ as search space.

Precisely, during the recycling-\IDRO-cycles \SRIDR constructs a solution in the search space
\begin{align}
\cU = \cK^\star_{J^\star}(\bA;[\bU,\br])\,, \label{eqn:Block_Recycling}
\end{align}
where $\br$ is the residual to the current \rhs and $\bU$ is the preimage of the reused auxiliary-vector-matrix $\bV$.
In the next chapter we will get to know a recycling principle that distinguishes from the one of this chapter by reusing the original search space $\cK_n(\bA;\bb\hia)$ for the construction of the solutions of the following systems.

\subsection{Termination of \SRIDR}

In this subsection we want to give an example in which one can observe the finite termination of \SRIDR. We underline that the examples given here shall not demonstrate practical usefullness but instead evidence that \IDR, \SRIDR and the method of Miltenberger \cite{IDR-Milten}, in the following called \textit{MI09}, in theory terminate within a finite number of iterations.

To realize a numerical behaviour that reflects these theoretrical properties we need to choose a simple system where the presence of numerical round-off does not ruin the theoretical finiteness of the methods. Thus we have chosen a small system dimension $N=100$. We chose $\bA = \operatorname{tridiag}(3,2,-1)$, $\bb_1 = \boldsymbol{1}$ and $\bb_2$ randomly from a normal distribution to ensure no correlation between $\bb_2$ and $\bb_1$.
\largeparbreak
We now solve the system for $\bb_1$ with \textsf{IDR}($20$) (i.e. $s=20$) - given as blue curves in figures \ref{fig:SRIDR_Finite1} and \ref{fig:SRIDR_Finite2} - and extract recycling data during this solution process at two different \IDRO-cycles $J^\star$ of the process - marked in the figures by \texttt{*}. Then we use \SRIDR and \textsf{MI09} with $s=20$ with the obtained recycling data to solve for $\bb_1$ and for $\bb_2$; so we solve for one \rhs where the recycling data is somehow correlated to and then for another \rhs which is completely arbitrary to this data.

\subsubsection{Termination Behaviour for $J^\star$ maximal} As $N=100$ and $s=20$ we know that \IDR terminates at least after $5$ \IDRO-cycles as $\cG_5 = \lbrace \bO \rbrace$. The same theoretically should hold for the recycling-\IDRO-cycles of \SRIDR.

\begin{figure}
\centering
\includegraphics[width=0.8\linewidth]{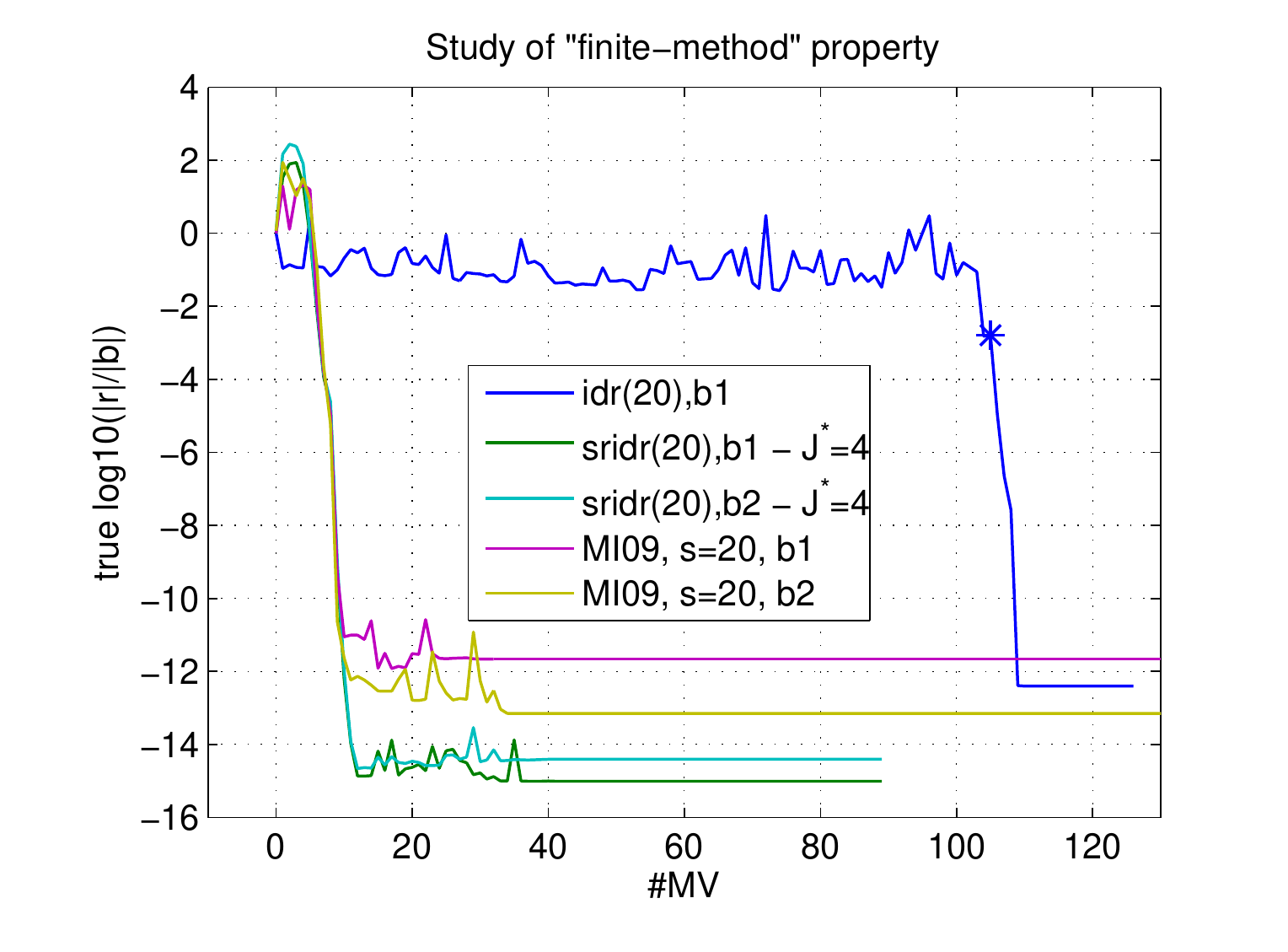}
\caption{Finite termination of \SRIDR (and seemingly \textsf{MI09}) after $5$ special \IDRO-cycles.}
\label{fig:SRIDR_Finite1}
\end{figure}

In figure \ref{fig:SRIDR_Finite1} we can see that \SRIDR indeed terminates nearly after $5$ recycling-\IDRO-cycles. This was predicted by the theory, no matter for which \rhs is solved.
We are surprised that in effect the method of Miltenberger behaves as good as \SRIDR although there does not exist any theory yet on how many iterations it needs to converge.

So we assumed that by coincidence it may chose the relaxation parameters near to those of \SRIDR. The relaxations can be found in table \ref{tab:RelaxComp1}. To judge the closeness of the relaxations one has to keep in mind that the order of the relaxations does not matter and that therefore even subsequent relaxations can suddenly restrict the residual into the highest old Sonneveld-space. However as \textsf{MI09} does not skip the auxiliary for-loop it always uses the same relaxtion for $s+1=21$ consecutive \MV-s. We have not analysed yet why this works.
\begin{table}
	\centering
\begin{tabular}{c|c|c}
 	index $j$ & \SRIDR & \textsf{Mi09} \\
	\hline
	\hline
	$\omega_1$ & 0.1877 & 0.2099 \\
	$\omega_2$ & 0.1844 & 0.2128 \\
	$\omega_3$ & 0.2057 & 0.2107 \\
	$\omega_4$ & 0.2055 & 0.2083 \\
	$\omega_5$ & 0.2372 & 0.2319 \\
	$\omega_6$ & 0.2051 & 0.2170 \\
	$\omega_7$ & 0.2133 & 0.2884
\end{tabular}
\caption{Comparison of relaxation coefficients for the convergence study from fig. \ref{fig:SRIDR_Finite1} for $\bb_1$.} \label{tab:RelaxComp1}
\end{table}

\subsubsection{Termination Behaviour for $J^\star = 3$}
\begin{figure}
\centering
\includegraphics[width=0.8\linewidth]{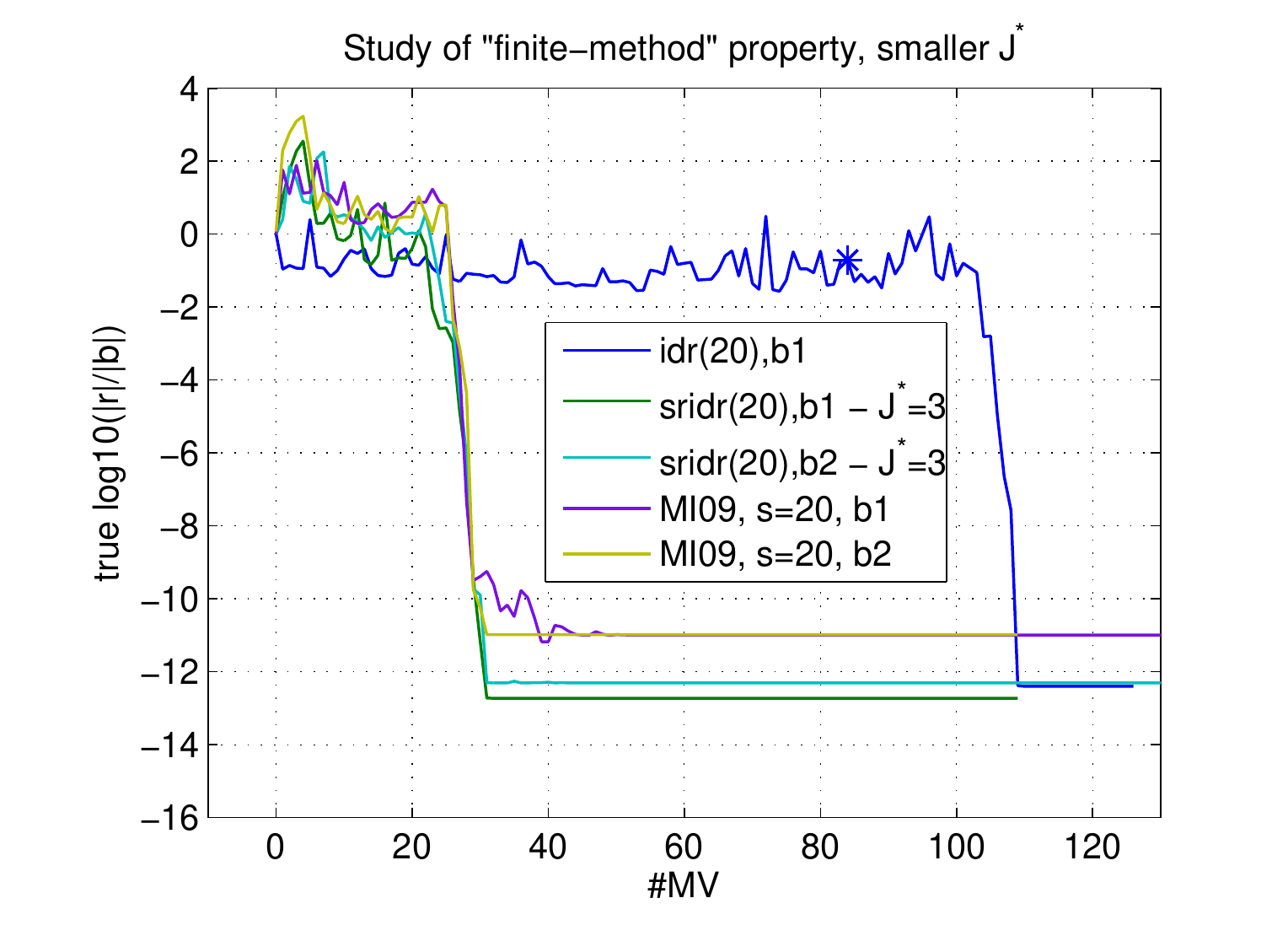}
\caption{Finite termination of \SRIDR (and seemingly \textsf{MI09}) for recycling data of lower Sonneveld-level.}
\label{fig:SRIDR_Finite2}
\end{figure}

In a second experiment we fetch the recycling data earlier after $J^\star = 3$ \IDRO-cycles of \textsf{IDR}($20$). One can see in fig. \ref{fig:SRIDR_Finite2} that the \texttt{*} has moved to the left by $21$ \MV-s. Thus in theory \SRIDR should now need these $21$ \MV-s additionaly compared to before. The experimental result confirms this. Also in this case, \textsf{MI09} turns out to converge surprisingly similar to \SRIDR.

\paragraph{Remark: Restrictions on $\boldsymbol{\omega}$ are necessary}
The above examples surprised us as for these and the used implementations the method \textsf{MI09} seems to terminate in the same step as \SRIDR. We do not know if there exist special implementations that allow alternative rules for choosing $\omega$ residual-dependent, such that \textsf{MI09} yields similar theoretical properties to \SRIDR.

\begin{figure}
\centering
\includegraphics[width=1\linewidth]{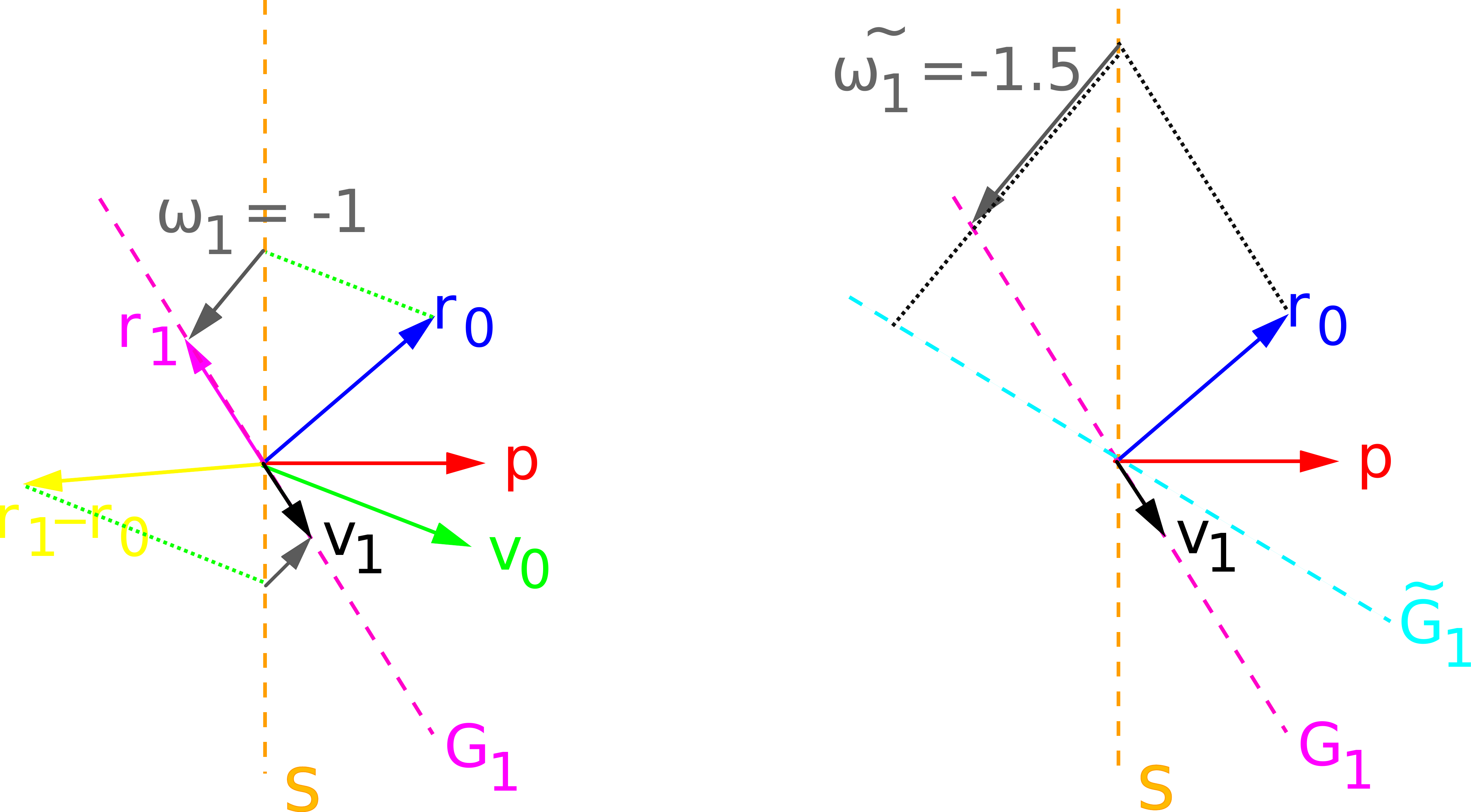}
\caption{Left: Conventional \IDReins; initial residual $\br_0$ and auxiliary vector $\bv_0$ in residual space. We arbitrarily chose $\omega_1 =-1$. With this $\br_1,\bv_1 \in \cG_1$ are computed. Right: Recycling of $\bv_1$ for same $\br_0$ with a different choice of the relaxation to $\tilde{\omega}_1 = -1.5$.}
\label{fig:Mi09_fail}
\end{figure}

However for general choice of $\omega$ we can give a simple counter-example that shows that in general \textsf{MI09} does not terminate whenever \SRIDR does. The example is given in fig. \ref{fig:Mi09_fail}.

Left in the figure we solve a $\R^2$-problem with initial residual $\br_0 \in \cG_0$. We first solve this problem with usual \IDReins without recycling. So we need one auxiliary vector $\bv_0 \in \cG_0$ and one shadow vector $\bp$, which defines $\cS$. In the solution process the residual is first orthgonalized with $\bv_0$ onto $\bp$ (green dotted line). Then the matrix-vector-product of this orthogonalized residual is computed (grey vector) and subtracted from the residual with $\omega_1 = -1$. By this we get the new residual $\br_1 \in\cG_1$ (magenta). Afterwards the next auxiliary vector is computed by
\begin{align*}
	\bv_1 = (\bI - \omega_1 \cdot \bA)\cdot(\bI-\bv_0 \cdot (\bp^H \cdot \bv_0)^{-1} \cdot \bp^H) \cdot (\br_1-\br_0)\,.
\end{align*}
The vectors $(\br_1-\br_0)$ in yellow, $-\bv_0 \cdot (\bp^H \cdot \bv_0)^{-1} \cdot \bp^H \cdot (\br_1-\br_0)$ green dotted, and $\bA \cdot(\bI-\bv_0 \cdot (\bp^H \cdot \bv_0)^{-1} \cdot \bp^H) \cdot (\br_1-\br_0)$ in grey, are plotted in the figure. This was just one cycle of the usual \IDReins-method. One can see that in the next orthogonalization of the residual on $\bp$ with $\bv_1$ we obtain $\br = \bO$, as both $\br,\bv$ are elements in $\cG_1$ and $\cG_1 \cap \cS = \lbrace \bO \rbrace$.

Now in the right part of the figure the same system is solved anew, this time with recycling. For this we hold $\bv_1,\bp$ as given from the left, as both \SRIDR and \textsf{MI09} use it. In the first step (for both \SRIDR and \textsf{Mi09}) $\br_0$ is orthogonalized w.r.t. $\bp$ by $\bv_1$ (grey dotted line from $\br_0$ to $\cS$). Then for this orthgonalized residual the matrix-vector-product is computed (resulting in the grey vector). This vector is scaled and reduced from the residual. In case of $\SRIDR$ we would choose this relaxation identical to $\omega_1$, as we would like to recycle the Sonneveld-space $\cG_1$. In \textsf{Mi09} instead, there is no restriction on the choice of the relaxiation. Intentionally, to get a differnce between \SRIDR and \textsf{Mi09}, for the latter method we choose $\tilde{\omega}_1 =-1.5$. By this the residual of \textsf{Mi09} is restricted into $\tilde{\cG}_1$. In difference to \SRIDR, we do not obtain $\br = \bO$ in the next iteration because $\br_1$,$\bv_1$ are not in the same Sonneveld-space.

\section{Short Recycling with Short Representations} \label{chap:ShortRepChap}
\SRIDR uses $J$ \MV-s and $3 \cdot s$ stored columns for the construction of a recycling solution with a residual which is orthogonal on $n = s \cdot J$ dimensions.
In this chapter instead a principle is presented, that needs $2 \cdot J$ \MV-s and $2 \cdot s$ stored columns for construction of a recycling solution with a residual which is orthogonal on $n = s \cdot J$ dimensions. In return it enables the recycling of the original search space $\cU = \cK_n(\bA;\bb\hia)$. Furthermore the principle from this chapter can be applied very simply to any type of Hessenberg-decomposition (e.g. \CG, \MINRES, \BiCG, \GMRES,...).

\subsection{Basic Idea: Block-Krylov-Matrices}
All Krylov-subspace methods that are known so far base on (generalized) Hessenberg-decompositions of the form
\begin{align*}
\bA \cdot \bV = \obV \cdot \obH\,,
\end{align*}
with $\bV \in \C^{N \times n}$. For large $n$ with $n = k \cdot J$ the basis $\bV$ cannot be stored completely; $\bH \in \C^{n \times n}$ and $\tbV = \bV(:,1:J:n) \in \C^{N \times k}$, depending on the choice of $k$, instead can be stored very well.

The idea is now to store merely the matrices $\tbV$ and $\bH$. This data is sufficient to reuse all information from the whole Hessenberg decomposition and to compute matrix-vector-products with $\bV$ and $\bV^H$ although $\bV$ is not stored explicitly. This is achieved through \textit{Short Representations}.
\begin{Theorem}[Short Representation] \label{thm:ShortReps}
	Let
	\begin{align*}
		\bA \cdot \bV = \obV \cdot \obH\,
	\end{align*}
	be a Hessenberg-decomposition with $n,J,k \in \N$, such that $k \cdot J = n$, with $\bH,\obH,\bV,\obV$ as in the notations, and $\tbV = \bV(:,1:J:n) \equiv [\bv_1,\bv_{1+J},\bv_{1+2\cdot J},...,\bv_{1+(k-1)\cdot J}] \in \C^{N \times k}$.

	Then one can compute and store a permutation matrix $\boldsymbol{\Pi} \in \C^{n \times n}$ in $\cO(n)$ and an upper right triangular matrix $\bK \in \C^{n \times n}$ in $\cO(n^3)$ (and in $\cO(n^2)$ if $\bH$ is banded), such that
	\begin{align*}
		\bV \cdot \bK &= [\tbV,\bA \cdot \tbV,...,\bA^{J-1} \cdot \tbV] \cdot \boldsymbol{\Pi}
	\end{align*}
	holds.
\end{Theorem}
\paragraph{Proof}
Construct $\bK$ with the following scheme:
\begin{align*}
\bK(:,i) &= \be_i,&\quad i &\in 1+J \cdot \N_0 \cap \lbrace 1,2,...,n\rbrace\\
\bK(:,i) &= \bH \cdot \bK(:,i-1),&\quad i &\in \lbrace 1,2,...,n\rbrace \setminus (1+J \cdot \N_0)
\end{align*}
Or, equivalently, as algorithm:
\begin{algorithmic}[1]
	\State $\bK := []$
	\For{$i=1,...,k$}
		\State $\bp := \be_{1 + (i-1) \cdot J} \in \C^n$
		\For{$j=1,...,J-1$}
			\State $\bK := [\bK,\bp]$
			\State $\bp := \bH \cdot \bp$
		\EndFor
		\State $\bK := [\bK,\bp]$
	\EndFor
\end{algorithmic}
Construct the permutation matrix $\bPi$, such that
\begin{align*}
	\bPi \cdot \be_{1 + i \cdot J + j} = \be_{1 + j \cdot k + i},\quad j=0,...,J-1,\ i=0,...,k-1\,. 
\end{align*}
The proposition now follows by simple calculations. q.e.d. A more detailed proof will be given in \cite{SRMR}.
\largeparbreak
\paragraph{Remarks}
\begin{itemize}
	\item We call $\hbV := [\tbV,\bA \cdot \tbV,...,\bA^{J-1} \cdot \tbV]$ \textit{block-Krylov-matrix}.
	\item Use of theorem \ref{thm:ShortReps}: Matrix-vector-products with $\hbV$ can be computed cheap and - instead of $\bV$ - usually $\tbV$ can be stored.
\end{itemize}

\begin{figure}
\centering
\includegraphics[width=1\linewidth]{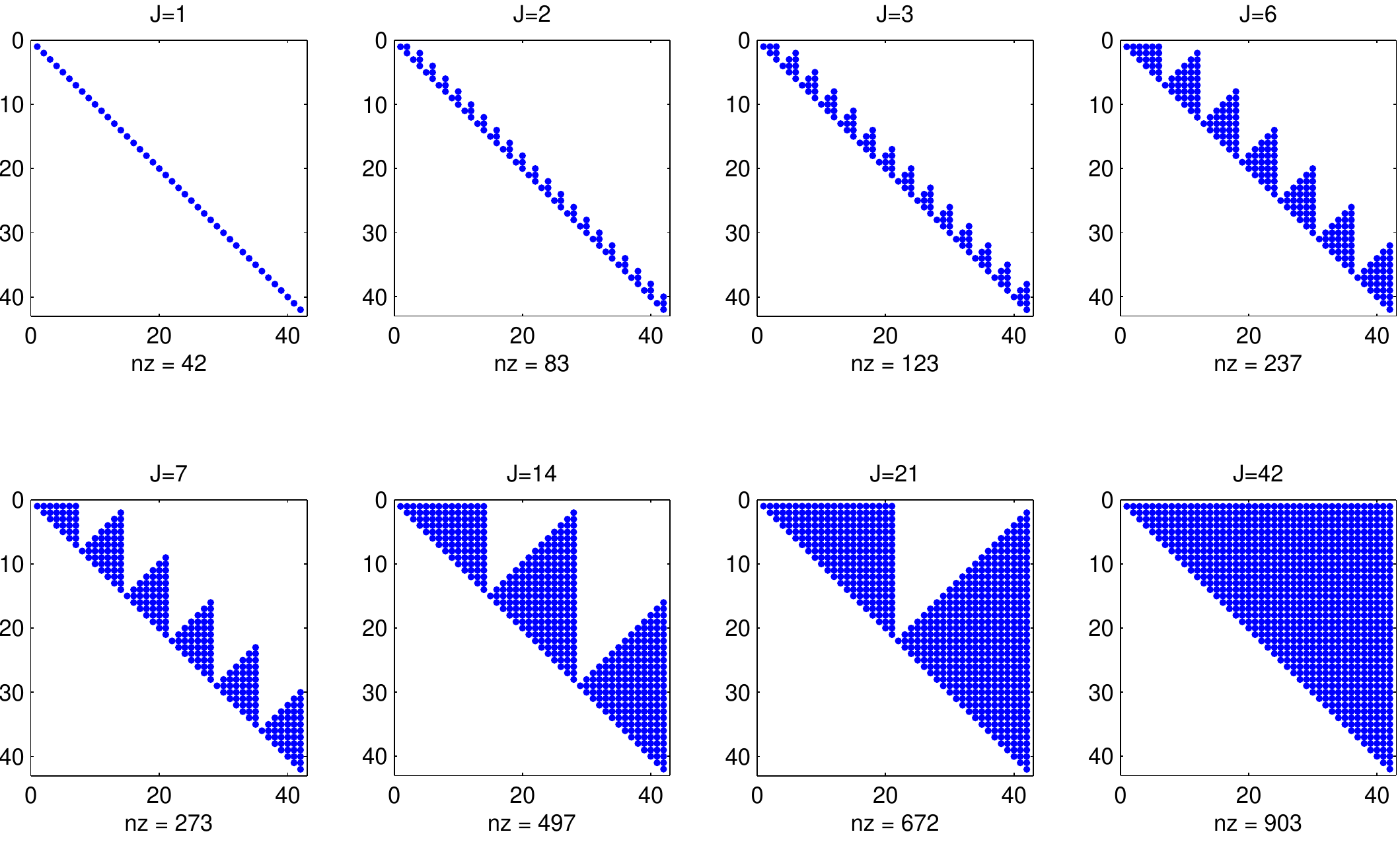}
\caption{Structure of short representation matrix $\bK$ for various choices of $J$.}
\label{fig:K_structure}
\end{figure}

To get an idea of how the matrix $\bK$ looks like, we exemplary computed $\bK$ for a tridiagonal matrix $\bT \in \R^{42 \times 42}$ ($n=42$) from a Lanczos decompostion for different choices of $J$. One can see in figure \ref{fig:K_structure} that due to the band structure of $\bT$ the matrix $\bK$ is also banded (with bandwidth $J$). As a remark, for tridiagonal $\bT$ the matrix $\bK$ can be computed in $\cO(n \cdot J)$.

From the band structure of $\bK$ one can deduce which columns $\bv_i$, $i=1,...,n$ from $\bV$ must be stored; namely those where only one element per column of $\bK$ is nonzero. For this the $J-1$ right-on columns of $\bK$ give the orthogonal relations of its power vectors $\bA^j \cdot \bv_i$ with the columns from $\bV$. For $J=1$ all columns of $\bV$ are stored, for $J=42$ only the first column of $\bV$ is stored.

\begin{figure}
\centering
\includegraphics[width=1\linewidth]{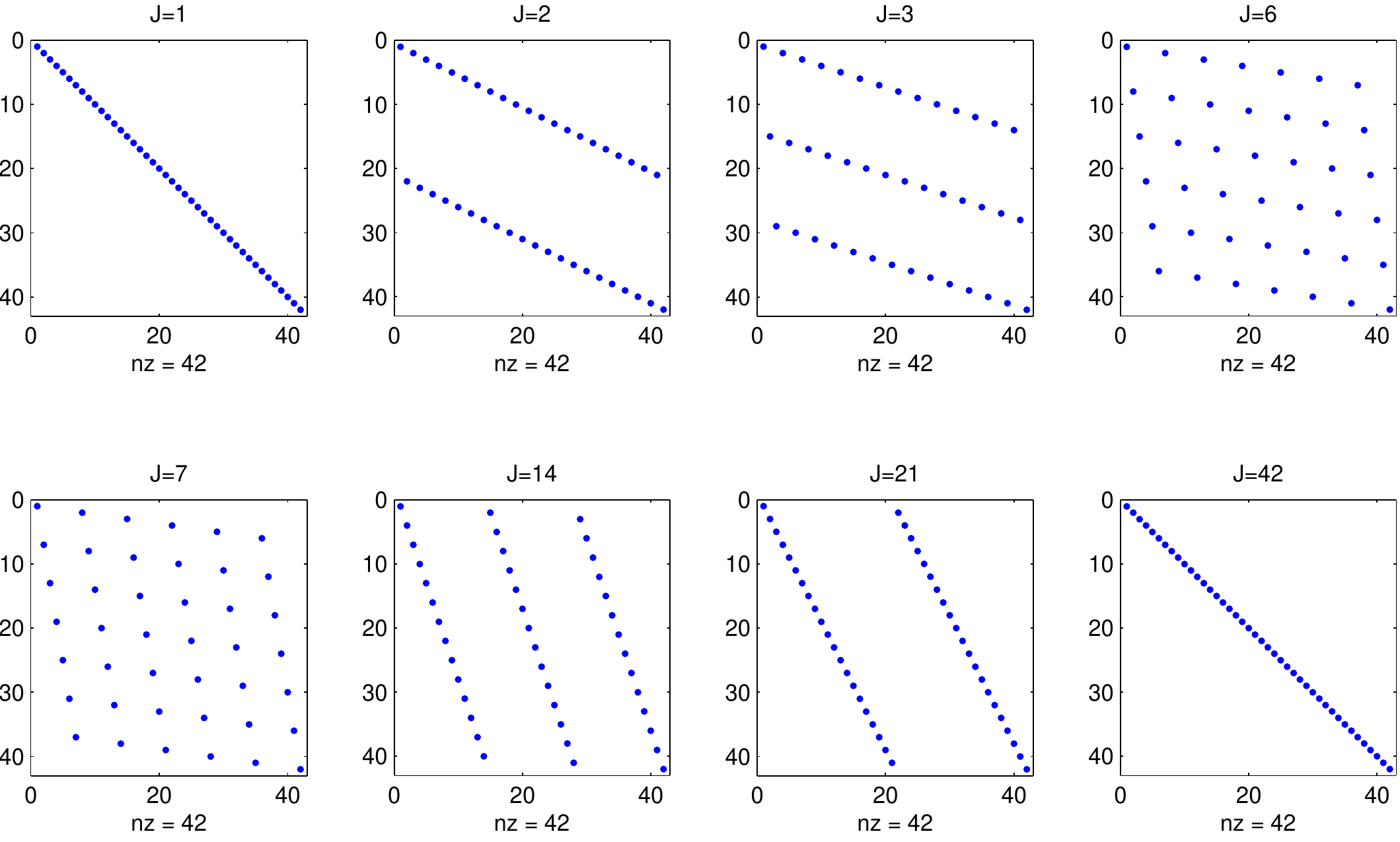}
\caption{Structure of permutation matrices $\bPi$ according to triangular matrices $\bK$ from fig. \ref{fig:K_structure} for various choices of $J$.}
\label{fig:P_structure}
\end{figure}

The structures of the permutation matrices $\bPi$ which belong to each triangular matrix $\bK$ from fig. \ref{fig:K_structure} are given in fig. \ref{fig:P_structure}. One can see some symmetry dependence on $J$, precisely $\bPi_{1/J} = \bPi^H_J$, where the sub-index indicates for which $J$ the matrix is constructed.

\subsection{Products with the Block-Krylov-Matrix} \label{chap:ProdsWithBKMs}

Matrix-vector-products with $\bV$ resp. $\bV^H$ can be computed in $J$ \MV-s with $\bA$ resp. $\bA^H$ and $J$ matrix-vector-products with $\tbV$ resp. $\tbV^H$. This can be done as follows:
One uses
\begin{align*}
	\bz \equiv \bV \cdot \by :=& \hbV \cdot \tilde{\by} \\
\tilde{\by} :=&\bPi \cdot \bK^{-1} \cdot \by
\end{align*}
for computation of $\bz$ from $\by$. The product with $\hbV$ is then computed via the Horner-scheme for $\tilde{\by} = (\tilde{\by}_1^H,\tilde{\by}_2^H,...,\tilde{\by}_J^H)^H$ by means of the following procedure:
\begin{algorithmic}[1]
	\State $\bz := \bO$
	\For{$j=0,...,J-2$}
	\State $\bz := \bz + \tbV \cdot \tilde{\by}_{J-j}$
	\State $\bz := \bA \cdot \bz$
	\EndFor
	\State $\bz := \bz + \tbV \cdot \tilde{\by}_{1}$
\end{algorithmic}

For products with $\bV^H$ one proceeds in the same way: First one reduces the problem to computing a matrix-vector-product with $\hbV^H$. This matrix-vector-product $\tilde{\by} = (\tilde{\by}_1^H,\tilde{\by}_2^H,...,\tilde{\by}_J^H)^H = \hbV^H \cdot \bz$ is then computed by means of the following scheme:
\begin{algorithmic}[1]
	\State $\bp := \bz$
	\For{$j=1,...,J-1$}
		\State $\tilde{\by}_j := \tbV^H \cdot \bp$
		\State $\bp := \bA^H \cdot \bp$
	\EndFor
	\State $\tilde{\by}_J := \tbV^H \cdot \bp$
\end{algorithmic}

\subsection{Use of Short Representations for Recycling} \label{chap:ShortRep_RecyclingMethods}
We now look at how short representations can be used to solve two consecutive systems.
We start by solving a first system
\begin{align*}
\bA \cdot \bx\hia &= \bb\hia\,.
\end{align*}
Using the data, that is computed in the solution process of the first system anyway, we try to compute an approximate solution to subsequent systems
\begin{align*}
\bA \cdot \bx\himu &= \bb\himu,\quad \mu = 1,2,...\,.
\end{align*}
The solutions to these systems that are obtained by using the recycling data are henceforth called \textit{recycling solutions}.

Many Krylov-subspace methods are based on a Hessenberg-decomposition of the form
\begin{align*}
	\bA\cdot \bV &= \obV \cdot \obH_1 \\
	\bB \cdot \bW &= \obW \cdot \obH_2\,.
\end{align*}
In this $\bV$ is a basis for the search space $\cU$ or the image space $\cV$, and $\bW$ is a basis for the test space $\cP$. For $\bV$ and $\bW$ one can construct short representations by virtue of the respective Hessenberg equation and given Hessenberg matrix.
\paragraph{Examples}
\begin{Aufz}
	\item For the (bi-)Lanczos procedure one has $\bH_2 = \bH_1^H$ and $\bB = \bA^H$.
	\item For the Arnoldi procedure it holds $\bB = \bA$, $\bH_2 = \bH_1$ and $\bW = \bV$.
\end{Aufz}
\paragraph{Remark}
There exist Krylov-subspace methods that keep a matrix $\bU \in \C^{N \times n}$, such that $\bA \cdot \bU = \bV$ holds. In this case $\bU$ is a basis of the search space $\cU$ and $\bV$ is a basis of the image space. Also for $\bU$ a short representation can be constructed whenever $\obH_1$ is known, as for $\bU$ the same Hessenberg equation holds as for $\bV$. Whenever a matrix $\bU$ is used, then $\bu_1 \parallel \bb\hia$ and $\cU = \opImage(\bU)$. If instead no $\bU$ is used, then it holds $\bv_1 \parallel \bb\hia$ and $\cU = \opImage(\bV)$.

The first approach where a matrix $\bU$ is used as basis of $\cU$ is called \textit{U-approach}, the second approach where instead $\bV$ is used as basis for $\cU$ is called \textit{V-approach}. (In general $\bU = \bV$  does not hold for $\bU$ from U-approach and $\bV$ from V-approach, as orthogonality of $\bV$'s columns differs from $\bU$'s.)
\largeparbreak
In the following I show how short representations can be used to recycle the information from a Hessenberg decomposition for the solution of a subsequent \rhs.

For this we go through the common Krylov-subspace methods and show how for each of them one can obtain a short recycling method by use of short representations. The short recycling methods are then called by the name of their parent Krylov-subspace methods with the prefix \textit{SR}, e.g. \SRBiCG for a short recycling method obtained from \BiCG.

\paragraph{Conjugate Gradients} Here the Hessenberg decomposition is the symmetric Lanczos decomposition (write $\bT := \bH$), in which the second equation is identical to the first one, thus $\bW = \bV$.
For \CG we use the V-approach. By means of short representations of $\bV$ recycling solutions for later \rhs-es $\bb\himu$, $\mu=1,2,...$ with \SRCG can then be computed by
\begin{align*}
	\bx\himu = \bV \cdot \big(\bT^{-1} \cdot (\bV^H \cdot \bb\himu)\big)
\end{align*}
Then it holds: $\br\himu \perp \opImage(\bV)$.
\paragraph{Minimal Residual} Here the Hessenberg decomposition is also the symmetric Lanczos decompositon und again $\bW = \bV$. There are two options now:
\begin{Aufz}
	\item In the V-approach one uses a short representation of $\bV$. Then a recycling solution for subsequent \rhs-es can be computed with
	 \begin{align*}
	 	\bx\himu 	&= \bV \cdot \obT^{\dagger} \cdot \obV^H \cdot \bb\himu\,,
	 \end{align*}
	 where the last column can be stored e.g. seperately. This approach has a disadvantage because \MINRES internally does not use an orthogonal basis $\bV$ of the search space but instead one of its image. This means that one cannot use this approach as a simultaneous computation process for a short representation and solution process of a first system. (In fact \MINRES internally uses the U-approach.) Due to this I propose a second option:
	 \item In the U-approach we instead store a short representation for $\bU$. The recycling solution can then be found by the last expression of
	 \begin{align*}
	 	\bx\himu 	&= \bU \cdot (\bV^H \cdot \bb\himu)\\
	 				&= \bU \cdot \bU^H \cdot \bA \cdot \bb\himu\,.
	 \end{align*}
\end{Aufz}
In both cases the recycling solution is residual optimal in the recycled original Krylov-subspace $\cK_n(\bA;\bb\hia)$ of the first \rhs. It holds: $\br\himu \perp \bA \cdot \cK_n(\bA;\bb\hia)$.
\largeparbreak
Both methods described so far (\SRCG, \SRMR) need storage for one $N \times k$-matrix and have computational cost of $2 \cdot J$ \MV-s with $\bA$ and $J$ matrix-vector-products with the stored $N \times k$-matrix and its transpose, respectively.

\paragraph{BiConjugate Gradients/Residual}
For this we only consider the case $\bA^H \neq \bA$. The Hessenberg decomposition is then the bi-Lanczos decomposition with in general $\bW \neq \bV$, and $\bH_2 = \bH_1^H$, $\bT := \bH_1$. As for \SRMR there exist two options:
\begin{Aufz}
	\item V-approach: For $\bV$ we store a short representation. The recycling solution is then
	\begin{align*}
		\bx\himu 	&= \bV \cdot \big(\bT^{-1} \cdot (\bW^H \cdot \bb\himu)\big)\,.
	\end{align*}
	\item U-approach: With a short representation for $\bU$ we obtain the recycling solution by
	\begin{align}
		\bx\himu 	&= \bU \cdot (\bW^H \cdot \bb\himu)\,. \label{eqn:RecSol_BiCG}
	\end{align}
\end{Aufz}
In both cases for the residual of the recycling solutions holds: $\br\himu \perp \opImage(\bW)$.
\largeparbreak
In contrast to the methods \SRCG and \SRMR above, for \SRBiCG one needs storage for \textit{two} instead of one $N \times k$-matrix. The costs instead for computing a recycling solution are the same as for \SRCG and \SRMR: $2 \cdot J$ \MV-s with $\bA$ and each $J$ matrix-vector-products with the stored $N \times k$-matrices $\tbU$ resp. $\tbV \in \C^{N \times k}$ (depending on U- or V-approach) and $\tbW^H \in \C^{k \times N}$.

Just as for \SRCG and \SRMR we also do not need \MV-s with $\bA^H$ for \SRBiCG, as we do only need to compute matrix-vector-products with the transpose of $\bW$.
\paragraph{Remark: Dual Systems}
The recycling data from \SRBiCG could also be used to find recycling solutions to dual systems: For the V-approach use
\begin{align*}
\bx\himu 	&= \bW \cdot \big(\bT^{-H} \cdot (\bV^H \cdot \bb\himu)\big)\,,
\end{align*}
for the U-approach instead use
\begin{align*}
\bx\himu 	&= \bW \cdot (\bU^H \cdot \bb\himu)\,.
\end{align*}
In the notation of the U-approach, for the residual of the recycling solution we then have: $\br\himu \perp \opImage(\bU)$.
\largeparbreak
Note that all recycling methods, that where mentioned in chapter \ref{chap:ShortRep_RecyclingMethods}, fullfill the condition: \rRD$ = k/2 \, \cdot \,$\rMV.

\paragraph{Remark: Efficiency}
One notices that the methods based on short representations need $2 \cdot J$ \MV-s with $\bA$, whereas \SRIDR for $s=k$ only needs $1 \cdot J$ \MV-s with $\bA$ to construct a recycling solution with a residual which is orthogonal on a $s \cdot J$-dimensional test space. Thus from a theoretical point of view \SRIDR seems to be more efficient.

I believe that the efficiency advantage of \SRIDR can be somehow translated into short recycling methods based on short representations.

\paragraph{\GMRES resp. \FOM}
Of course one could also use short representations for recycling the information of an Arnoldi decomposition. The formulas for the computation of the recycling solutions are then just like those of \SRMR resp. \SRCG. However for large $n \gg k$, where $k$ is the number of storable columns, one cannot compute the Hessenbergmatrix $\obH$ in $\cO(n)$ \MV-s with $\bA$.

\paragraph{Remark: Exceptional Recources}
If a short recycling method is used for example on a cluster computer, then one could think of the approach to compute the Hessenberg matrix once with large memory consumption, then store the short representations, and release the used memory space afterwards.
Then $\bH$ could be used e.g. for a \SRGMRES method.

\section{Increasing Robustness for Short Representations}\label{chap:Robust}
\subsection{Stability Issues}
From numerical experiments one notices that the short representations are sensitive to round-off errors. Of course this is intrinsic, as we use a structural compression strategy (based on the structure of the Krylov-subspaces), but one can distinguish two error sources for the recycling solutions.
\begin{Aufz}
	\item For increasing dimension $n$ of the recycled bases the numerically computed matrices $\bV$,$\bW$ loose their (bi-) orthonormality, thus the property $\bW^H \cdot \bV=\bI$ resp. $\bV^H \cdot \bV = \bI$ gets lost.\\
	When these matrices are now used for computation of recycling solutions, then the recycling solutions are already unprecise due to these inaccuracies. Henceforth this influence is called \textit{data-accuracy}.
	\item Apart from the recycled basis matrices there occur round-off errors due to the short representations which are computed to store the bases. Henceforth we call this influence \textit{representation-accuracy}.
\end{Aufz}
In reality the effects of both bad data and representation accuracy overlap.

The impact of bad data-accuracy can be handled by iterative orthogonalization of the residual on small blocks of columns for which the (bi-)orthogonality nearly holds. The representation accuracy on the other hand increases, when $J$ is chosen smaller. This is obvisiously, as in the limit for $J=1$ the representation accuracy is optimal since for $J=1$ there exists no short representation. However the representation-accuracy also increases when the condition of the Hessenberg matrix gets smaller, what always happens for smaller choice of $n$.

\largeparbreak
\subsection{Basic Idea: Stabilization by Blocking}
The idea consists of dividing a basis matrix $\bV$ into smaller blocks $\bV_1,\bV_2,...,\bV_\ell$, $\ell \in \N$, and storing short representations for each block.

As the blocks are smaller, it can be expected that the data-accuracy of the respective blocks increases (i.e. better orthogonality of $\bV^{(i)}$ resp. biorthogonality of $\bW^{(i)},\bV^{(i)} \in \C^{N \times n_i}$, $i=1,...,\ell$). On the other hand one can also expect, that due to smaller dimensions $n_i$ of each block the representation accuracy for each block increases. (Summing up one cures two problems in one approach.)
\largeparbreak
In the following I exemplary show for \SRBiCG, how this partitioning in blocks works in practice.

From the solution process of one system with \BiCG and $n$ \MV-s with both $\bA$ and $\bA^H$ one computes columns of matrices $\bU,\bV,\bW \in \C^{N \times n}$, such that
\begin{align*}
	\bA \cdot \bU &= \bV\\
	\bA \cdot \bV &= \obV \cdot \obT\\
	\bA^H \cdot \bW &= \obW \cdot \ubT^H\\
	\bW^H \cdot \bV &= \bI\,.
\end{align*}
Consider the subdivision of $\bT$ in $\bT = \operatorname{diag}(\bT^{(1)},\bT^{(2)},,...,\bT^{(\ell)}) + \textbf{E}$ (where the $\bT^{(i)}$ are just the diagonal blocks of $\bT$ and $\textbf{E}$ is the rest), and respective subdivisions of $\bU,\bV,\bW$ into $\bU = [\bU^{(1)},\bU^{(2)},...,\bU^{(\ell)}]$ with $\bU^{(i)} = [\bu^{(i)}_1,\bu^{(i)}_2,...,\bu^{(i)}_{n_i}]$, $\obU^{(i)} = [\bu^{(i)}_1,\bu^{(i)}_2,...,\bu^{(i)}_{n_i},\bu^{(i)}_{n_i+1}]$ and the same for $\bV^{(i)}$,$\bW^{(i)}$,$\obV^{(i)}$,$\obW^{(i)}$. To circumvent applications with $\textbf{E}$ we define the matrices
\begin{align*}
\bA_{\bU,i+1} &:= \big(\bI - \bu^{(i)}_{n_i} \cdot (\tilde{\bw}^{(i)}_{n_i})^H\big) \cdot \bA\\
\bA_{\bV,i+1} &:= \big(\bI - \bv^{(i)}_{n_i} \cdot (\bw^{(i)}_{n_i})^H\big) \cdot \bA\\
\bA_{\bW,i+1} &:= \bA \cdot \big(\bI - \bv^{(i)}_{n_i} \cdot ({\bw}^{(i)}_{n_i})^H\big)
\end{align*}
with $\bA_{\bU,1} = \bA_{\bV,1} = \bA_{\bW,1} = \bA$ and $\tilde{\bw}^{(i)}_{n_i}:=\bA^H \cdot \bw_{n_i}^{(i)}$. Then the following recursions hold:
\begin{align*}
\bA_{\bU,i} \cdot \bU^{(i)} &= \obU^{(i)} \cdot \obT^{(i)}\\
\bA_{\bV,i} \cdot \bV^{(i)} &= \obV^{(i)} \cdot \obT^{(i)}\\
\bA_{\bW,i}^H \cdot \bW^{(i)} &= \obW^{(i)} \cdot \big(\ubT^{(i)}\big)^H\,.
\end{align*}
Thus for each $\bU^{(i)}$,$\bV^{(i)}$,$\bW^{(i)}$, $i=1,...,\ell$ one can store short representations.

These representations can then be used to orthogonalize the residual iteratively on each $\opImage(\bW^{(i)})$, $i=1,...,\ell$. As an example we show this with the U-approach:
\begin{algorithmic}[1]
	\State \textit{//pre-conditions: $\br = \bb - \bA\cdot \bx$, $(\bW^{(i)})^H \cdot \bA \cdot \bU^{(i)} = \bI_{n_i}$}
	\For{$i=1,2,...,\ell$}
		\State $\bc := (\bW^{(i)})^H \cdot \br$\quad \textit{//internally computed via short rep. of $\bW^{(i)}$}
		\State $\bx := \bx + \bU^{(i)} \cdot \bc$\quad \textit{//internally computed via short rep. of $\bU^{(i)}$}
		\State $\br := \bb - \bA \cdot \bx$
	\EndFor
	\State \textit{//post-conditions: $\bx \in \opImage(\bU) \,:\, \br \perp \opImage(\bW)$}
\end{algorithmic}
Theoretically, line 5 would not be necessary in case of perfect data-accuracy as then $\bW^H \cdot \bV = \bI$ would hold. A sufficiently small number of columns $n_i$ for $\bU^{(i)},\bW^{(i)} \in \C^{N \times n_i}$ ensures a good representation-accuracy (e.g. consider $n_i=1$).

\section{A-posteriori-Improvement for Recycling Solutions}\label{chap:ApostImprov}
So far I have described how information of a former system can be reused under small storage and computation cost. When a solution, e.g. $\bx\himu$ for some $\mu \in \N$ with \SRBiCG by \eqref{eqn:RecSol_BiCG}, has been constructed with a short recycling approach, then now the question arises, how this solution can be improved in sense of further reduction of the residual norm by further iterations.

A naive approach consists of using $\bx\himu$ as initial guess for a conventional Krylov-subspace method such as \GMRES or \MINRES. However with this one looses the orthogonality of $\bx\himu$'s residual $\br$ to the recycled test space: $\br \perp \cP$. We want to conserve this property, which is called \textit{a-posteriori-orthogonality} in the remainder of the text.
\largeparbreak
In this chapter I am going to present a principle for iterative improvement of the recycling solution, which enables a-posteriori-orthogonality. This principle leads to a method, that
\begin{Aufz}
	\item is based on \IDRO,
	\item works with short recurrences,
	\item does not need transpose-products (i.e. with $\bA^H$),
	\item and is efficient, i.e. by each \MV with $\bA$ both dimensions of the search and test space are nearly increased by one.
\end{Aufz}

\subsection{Basic Idea: \IDR-Cycles to conserve A-posteriori-Orthogonality}
Let us revisit \SRIDR: After computing $J^\star$ \IDRO-cycles with \SRIDR, the recycling solution with $\br \in \cG_{J^\star}$ is found. By simply continuing with usual \IDR this recycling solution can be improved without loosing orthogonality properties.

Precisely for the residual and the auxiliary vectors the property $\br,\bv_1,...,\bv_s \in \cG_{J^\star}$ is conserved.

From \eqref{eqn:Sonneveld_Gut} we can see, that a Sonneveld-space $\cG_J$ is just the orthogonal complement of a testspace on which the residual is orthogonalized. Thus the test spaces of \IDRO-methods are:
\begin{align*}
	\cP_{\text{IDR},J} = \Omega_J(\bA^H) \cdot \cK^\star_{J}(\bA^H;[\bp_1,...,\bp_{s}])\,.
\end{align*}
As a summary to a-posteriori-iterations for \IDRO-methods, these methods automatically ensure a-posteriori-orthogonality of the residual.
\largeparbreak
Now we come to methods based on short representations:
These methods, which were described in chapter \ref{chap:ShortRepChap}, do not orthogonalize the residual on $\cP_{\text{IDR}}$, but instead on  a (block-) Krylov-subspace of the form
\begin{align}
\cP_{\text{SR},J} = \cK^\star_J(\bA^H;[\bw_1,...,\bw_{k}])\,.  \label{eqn:BlockTestSpace2}
\end{align}
Actually we only considered the case $k=1$ so far.

The question arises, if the a-posteriori-orthogonality principle from \SRIDR can be transferred to these methods.
For this first of all we notice that $\cP_{\text{SR},J}$ and $\cP_{\text{IDR},J}$ are quite similar.
\largeparbreak
By using a method of chapter \ref{chap:ShortRepChap} one can realize, that at the end of the computation of a recycling solution, the following data $\br\himu,\bv_1,...,\bv_s \perp \cP_{\text{SR},J}$ exists for some $s \in \N$.

The basic idea is now to modify these vectors $\br\himu,\bv_1,...,\bv_s \perp \cP_{\text{SR},J}$, such that they are elements of $\cG_{J^\star}$, resp. orthogonal w.r.t. $\cP_{\text{IDR},J^\star}$ (under appropriate choices of $\Omega$, $J^\star$ and $[\bp_1,...,\bp_{s}]$). This is useful in the sense, that afterwards one can apply \IDRO-cycles for the modified residual and restrict it into subspaces of $\cG_{J^\star}$. As above for \SRIDR, this ensures conservation of $\br \perp \cP_{\text{IDR},J^\star}$. Of course $\dim(\cP_{\text{IDR},J^\star}) \approx \dim(\cP_{\text{SR},J})$ is desirable in order to loose not too many orthogonality properties.

\subsection{Transforming Block-Krylov-Subspaces into Sonneveld-Spaces} \label{chap:Trans_K_to_G}
In the following I present a couple of approaches for construction of appropriate data in an eligable Sonneveld-space. In this we consider the situation where the space $\cP_{\text{SR},J}$ is given by data $J,\bw_1,...,\bw_{k}$, and a space $\cP_{\text{IDR},J^\star}$ shall be defined by a useful choice of $J^\star,\bp_1,...,\bp_{s},\Omega_j(\cdot)$.

\paragraph{IDR(1) for \BiCG}
In the bi-Lanczos decomposition with left starting vector $\bw_1$ we obtain a testspace $\cP_{\text{SR},J}$ from \eqref{eqn:BlockTestSpace2} with $k=1$ and $J=n$.

In this decomposition it holds for the last column of $\obV$: $\bv_{n+1} \perp \opImage(\bW) \equiv \cP_{\text{SR},n}$.
We can choose the polynomial $\Omega_J$ of degree $J$ with non-zero leading coefficient arbitrarily and then compute $\hbw = \Omega_J(\bA^H) \cdot \bw_1$. I propose to choose this polynomial such that $\hbw = \bw_{n+1}$, where $\bw_{n+1}$ is the last column of $\obW$. (We actually do not need to know this polynomial explicitly in the sense of coefficients or roots).

Now assume that we have computed a recycling solution with \SRBiCG, then we have the residual $\br \perp \cP_{\text{SR},n}$. We can define $\bt := \Omega_J(\bA) \cdot \br \in \cG_J$, where $\cG_J$ is defined by $J=n,s = k = 1$ and $\bp_1 = \bw_1$. We also define $\bg := \Omega_J(\bA) \cdot \bv_{n+1} \in \cG_J$. Theoretically, for $\bt,\bg$ we could now do \IDRO-cycles by
\begin{align}
\begin{split}
\bt_{j+1} :=&(\bI - \omega_{j+1} \cdot \bA) \\
&\cdot \big(\bI -  \bg_j \cdot (\bw_1^H \cdot \bg_j)^{-1} \cdot \bw_1^H \big) \cdot \bt_j \fuer j=J^\star,J^\star+1,...
\end{split}\label{eqn:IDR_t}	\\
	\begin{split}
	\bg_{j+1} :=&(\bI - \omega_{j+1} \cdot \bA)  \\
	&\cdot \big(\bI -  \bg_j \cdot (\bw_1^H \cdot \bg_j)^{-1} \cdot \bw_1^H \big) \cdot (\bt_{j+1} - \bt_j)\fuer j=J^\star,J^\star+1,...\,. \label{eqn:IDR_g}
	\end{split}
\end{align}
In this we choose $\omega_{j+1} \in \C\setminus\lbrace 0 \rbrace$, $j = J^\star,...$ . For the objects computed with these expressions hold: $\bt_j,\bg_j \in \cG_j$, $j\geq J^\star$.

We can write two related expressions
\begin{align}
\begin{split}
\br_{j+1} :=&(\bI - \omega_{j+1} \cdot \bA) \\
&\cdot \big(\bI -  \bv_j \cdot (\hbw^H \cdot \bv_j)^{-1} \cdot \hbw^H \big) \cdot \br_j \fuer j=J^\star,J^\star+1,...
\end{split}\label{eqn:IDR_r}	\\
\begin{split}
\bv_{j+1} :=&(\bI - \omega_{j+1} \cdot \bA)  \\
&\cdot \big(\bI -  \bv_j \cdot (\hbw^H \cdot \bv_j)^{-1} \cdot \hbw^H \big) \cdot (\br_{j+1}-\br_j) \fuer j=J^\star,J^\star+1,...\,.
\end{split}\label{eqn:IDR_v}
\end{align}
The relation is as follows: When multiplying \eqref{eqn:IDR_r} from left with $\Omega_J(\bA)$, this yields \eqref{eqn:IDR_t}. The same holds for \eqref{eqn:IDR_v} and \eqref{eqn:IDR_g}. The scalar products $\hbw^H \cdot \bv_j$ and $\hbw^H \cdot \br_j$ are identical respectively to $\bw_1^H \cdot \bg_j$ and $\bw_1^H \cdot \bt_j$ for all $j\in J^\star + \N_0$.

By this relations, from the restriction of $\bt_{j+1} \in \cG_{j+1}$ of $s=1$ fewer dimensions follows an orthogonality of $\br_{j+1} \perp \cP_{\text{SR},j+1}$ on the testspace increased by $s=1$ dimension.

\largeparbreak
For practical computations we just use expressions \eqref{eqn:IDR_r} and \eqref{eqn:IDR_v}. By this we get a scheme to orthogonalize $\br$ by computation of two \MV-s with $\bA$ onto one additional test space dimension.

The reason why I presented this rather unintuitive approach was that we obtain $\dim(\cP_{\text{IDR},J^\star}) = \dim(\cP_{\text{SR},J})$ (precisely $\cP_{\text{IDR},J} = \Omega_J(\bA^H) \cdot \cP_{\text{SR},J}$, $J=J^\star=n$), which means no loss of orthogonality information, and we do not need to compute additional data for preparing the a-posteriori-iterations because we already have $\hbw = \bw_{n+1}$ from the bi-Lanczos decomposition.

\largeparbreak
To increase the efficiency of a-posteriori-iterations, we would like to replace \IDReins by \IDR with $s>1$, because for $s>1$ \IDR seems more efficient, cf. \IDRO-references.

\paragraph{\IDR for \BiCG}
To be able to choose $s>1$, we need additional vectors $\bp_2,...,\bp_s \in \C^N$ for construction of the Sonneveld-space, and additional auxiliary vectors to be able to compute \IDRO-steps. For the auxiliary vectors we can increase the number of computed basis vectors of $\bV$ from the bi-Lanczos decomposition, to obtain the columns $\bv_{n+2},...,\bv_{n+s}$. These columns can be used as additional auxiliary vectors because also $\bv_{n+2},...,\bv_{n+s} \perp \opImage(\bW) \equiv \cP_{\text{SR},J}$, where $J$ and $\cP_{\text{SR},J}$ are as above. Also as above, we choose $\bp_1 = \hbw$. This ensures conservation of the residuals orthogonality to $\cP_{\text{SR},J}$ for all following \IDRO-cycles. The other shadow vectors $\bp_2,...,\bp_s$ are chosen randomly. By analysing this construction one obtains that if we compute \IDRO-cycles with $\br,\bv_{n+1},...,\bv_{n+s}$ and shadow vectors $\bp_1,...,\bp_s$, the residual remains orthogonal to the recycled test space.

The reason why this works is that one could replace $\cS$ in the recursive definition of $\cG_j$ by $\cS_j \supset \cS_{j+1}$ and still have $\cG_j \subset \cG_{j+1}$, $j=0,...,J-1$. Then one would show that the Krylov-subspaces of the vectors $\bp_2,...,\bp_s$ by random choice satisfy certain conditions with almost 100\% probability. However these discussions are beyond the scope of this text.

\paragraph{Remark}
Later on we will see a method based on a multiple-left bi-Lanczos decomposition. There we will see that both the cases $s=1,k=1$ and $s>1,k=1$ are special cases of a more general scheme, suited for \IDR and \MLBiCG with $s \geq k \geq 1$.

\section{Short Recycling for \MLBiCG}

Until this chapter the following methods were presented: \SRIDR and \SRBiCG for non-Hermitian systems and \SRMR (\SRCG) for Hermitian (positive definite) systems.

In all of these methods we first solve one system
\begin{align*}
 	\bA \cdot \bx\hia = \bb\hia
\end{align*}
with a conventional method, e.g. \IDR, \BiCG, \MINRES or \CG. From this solution process we extract recycling data for construction of short recycling solutions to the following systems.

For all methods except \BiCG the relation \rRD $\approx$ \rMV holds. For \BiCG instead one needs \MV-s with $\bA$ and $\bA^H$ for both construction of the test space and the search space. Precisely for \BiCG it holds: \rRD$=\lceil 1/2 \cdot$\rMV$\rceil$. This is inefficient, compared e.g. to \IDR: \rRD$=s \cdot \lceil 1/(s+1) \cdot$\rMV$\rceil$.

Due to this inefficiency of \BiCG the methods \MLBiCG and \IDR were developed as remedy.

As for \IDR I already have presented a short recycling method in chapter \ref{chap:SonnRecChap}, we now deal with a short recycling method based on \MLBiCG, which we will call, in analogy to the others, \SRMLBiCG. In the following $k \in \N$ is user-defined, as $s \in \N$ is in \IDR.

Compared to \SRIDR, \SRMLBiCG has the advantage to recycle the original search space $\cK_n(\bA;\bb\hia)$. As a drawback, for the solution of the first system $(\iota)$ the method \SRMLBiCG needs slightly more \MV-s than \SRIDR (precisely $2 \cdot J$) to compute both the recycling data and a solution to the first system.

Before we go on we want to underline that the method which is presented in this chapter has not been implemented so far. There is no evidence for practical usefulness yet, i.e. it could turn out that the geometrical construction is too sensitive to numerical round-off.

\subsection{Recycling of a transpose-free multiple-left bi-Lanczos Decomposition}

\paragraph{Motivation}
In the following we want to build a scheme for computation of a multiple-left bi-Lanczos decomposition (\MLbiL). We will use the information of this decomposition for computation of recycling solutions for subsequent systems.

To deal with \MLbiL we clearify the following questions in this order:
\begin{itemize}
	\item What form does the multiple-left bi-Lanczos decomposition have?
	\item How can one use it to compute recycling solutions?
	\item Which part of the \MLbiL is actually needed for recycling computations?
	\item How can this part be computed efficiently?
\end{itemize}

\paragraph{Review (multiple-left bi-Lanczos)}
In the U-approach, the method \MLBiCG is based on the following form of a multiple-left bi-Lanczos decomposition
\begin{align*}
	\bA \cdot \bU &= \bV\\
	\bA \cdot \bV &= \obV \cdot \obT\\
	\bA^H \cdot \bW &= \obW \cdot \ubT^H\,
\end{align*}
with $\bW^H \cdot \bV = \bI$, in which $\bT$ is not a tridiagonal matrix but instead has $k$ upper subdiagonals. As usual the columns of $\bV$ build a basis of the Krylov-subspace $\cK_n(\bA;\bA \cdot \bu_1)$. However the columns of $\bW$ build a basis of the block-Krylov-subspace $\cK_n(\bA^H;[\bw_1,...,\bw_k])$ with $\opImage(\bW(:,1:i)) = \cK_i(\bA^H;[\bw_1,...,\bw_k])$, $\forall i=1,...,n$. At the beginning of the method $\bu_1$ and $\bw_1,...,\bw_k \in \C^N$ are chosen arbitrarily by the user. Further details can be found in \cite{BiL-Block}.

For simplicity we assume in the following, that $\rank(\bW) = n$, where $n$ is the size of the decomposition, and that in the bi-Lanczos procedure no break-down occurs.

\paragraph{Computation of (recycling) solutions with the multiple-left bi-Lanczos decomposition}
The expression for computation of recycling solutions in case of the multiple-left bi-Lanczos decomposition is identical to the one for \SRBiCG, cf. \eqref{eqn:RecSol_BiCG}, as the \MLbiL has the same form as the conventional bi-Lanczos decomposition and satisfies the same biorthogonality properties.

For $\bU$ and $\bW$ one can store short representations. We will see in later on how this can be done. By being able to compute products with $\bU$ und $\bW^H$, one can approximately solve linear equation systems  (here exemplary in the U-approach) as follows:
\begin{align*}
	\bA \cdot \bx\himu &= \bb\himu \\
	\Rightarrow \bW^H \cdot \bA \cdot \bx\himu &= \bW^H \cdot \bb\himu\\
	\Leftarrow \bW^H \cdot \bA \cdot \bU \cdot \by &= \bW^H \cdot \bb\himu\\
	\Leftrightarrow \bx\himu &= \bU \cdot \bW^H \cdot \bb\himu\,,
\end{align*}
thus identical to the expression \eqref{eqn:RecSol_BiCG} from \SRBiCG.

\paragraph{Remark: Calculation Procedure}
With the method presented in the following for computation of short representations for $\bU$ and $\bW$ I believe that it is not possible to simultaneously solve a first system $(\iota)$. Therefore the whole calculation procedure for the solution of a sequence is as follows:
\begin{Aufz}
	\item For the first system $\bA \cdot \bx\hia = \bb\hia$ we compute a \MLbiL for which we choose the first column $\bu_1$ of $\bU$ parallel to $\bb\hia$ (in the U-approach). The vectors $\bw_1,...,\bw_k \in \C^N$ are chosen randomly as for \IDR. A proper choice of $n$ (size of the \MLbiL) is estimated with the help of former systems.
	\item From the \MLbiL we extract and store recycling data for products with $\bU$ and $\bW^H$ (in the U-approach).
	\item The short representations are then used to solve the first system for $\bx\hia$.
	\item If $\bx\hia$ is not precise enough, then a few a-posteriori-iterations are performed (explained in a later chapter).
	\item Afterwards the following systems $(\iota+\mu)$, $\mu=1,2,...$ are solved just like the first one by doing steps 3 and 4.
\end{Aufz}

\subsection{Computation of a multiple-left bi-Lanczos decomposition}
I now present how the short representations of $\bV$ and $\bW \in \C^{N \times n}$ and the band matrix $\bT \in \C^{n \times n}$ from a \MLbiL of size $n=k \cdot J$ can be computed in $n + n/k$ \MV-s with $\bA$ and without using products with $\bA^H$. Obviously from the computation of a short representation for $\bV$ a short representation for $\bU$ can be recaptured.

\subsubsection{Basic Idea: Level-Vectors}
The test space, that is spanned by $\bW$, is
\begin{align*}
	\cP &= \operatorname{span}\lbrace \bw_1,...,\bw_k,\bA^H \cdot \bw_1,...,\bA^H \cdot \bw_s,...,(\bA^H)^{J-1} \cdot \bw_k \rbrace\\
		&= \operatornamewithlimits{span}_{i=1,...,n}\lbrace (\bA^H)^{g_k(i)} \cdot \bw_{r_k(i)}\rbrace
\end{align*}
with $g_k(i) := \lfloor (i-1)/k \rfloor$ and $r_k(i) := \big((i-1)\,\operatorname{mod}\,k+1\big) \leq k$, cf. \cite{BiL-TF-Block}.
\largeparbreak
I now explain how the \MLbiL can be computed. For this we define $\bp_i := (\bA^H)^{g_k(i)} \cdot \bw_{r_k(i)}$, $i=1,...,n$. To use the short recursion
\begin{align}
	\bv_{i+1} := \bA \cdot \bv_i - \sum_{c=\max(i-k,1)}^{i} \bv_c \cdot t_{i,c}\quad \perp \opspan\lbrace\bw_{\max(i-k,1)},...,\bw_i\rbrace \label{eqn:LinComb_Lancz}
\end{align}
in the $i$-th iteration of \MLbiL for computation of the next column of $\bV$, the vector $\bA \cdot \bv_i$ can also be orthogonalized w.r.t. $\bp_{\max(i-k,1)},...,\bp_i$ instead of $\bw_{\max(i-k,1)},...,\bw_i$.

As we cannot compute $\bp_i$ for $i>s$ due to avoidance of products with $\bA^H$, all scalar products of the form $\langle \bp_c\,,\,\bv_i\rangle$, that occur during the multiple-left bi-Lanczos procedure for computation of the linear factors $t_{c,i}$, must be replaced as follows:
\begin{align*}
	\langle \bp_c\,,\,\bv_i \rangle = \langle \bw_{r_k(c)}\,,\,\bA^{g_k(c)} \cdot \bv_i\rangle\,,\quad \forall i,c=1,...,n\,,
\end{align*}
cf. \cite{BiL-TF-Block}. As we do not want to compute up to $g_k(n)=J$ \MV-s with $\bA$ for each column $\bv_i$, $i=1,...,n$ of $\bV$, we use \textit{level-vectors} $\bv^{(j)}_i := \bA^{j-1} \cdot \bv_i$. This is comparable with the approach in \cite{MLBiCGstab}.

\subsubsection{Constructing a Short-Recurrence-Procedure}
The first $k$ columns $\bv_1,...,\bv_k$ of $\bV$ can be computed identically to the conventional bi-Lanczos procedure (but with recurrence length $k+1$ instead of $2$), as the columns $\bw_1,...,\bw_k$ are explicitly given by the user. We scale $\bv_1,...,\bv_k$ such that $\langle\bv_j,\bw_i\rangle = \delta_{i,j}$, $i,j=1,...,k$ holds.

\begin{figure}
\centering
\includegraphics[width=1\linewidth]{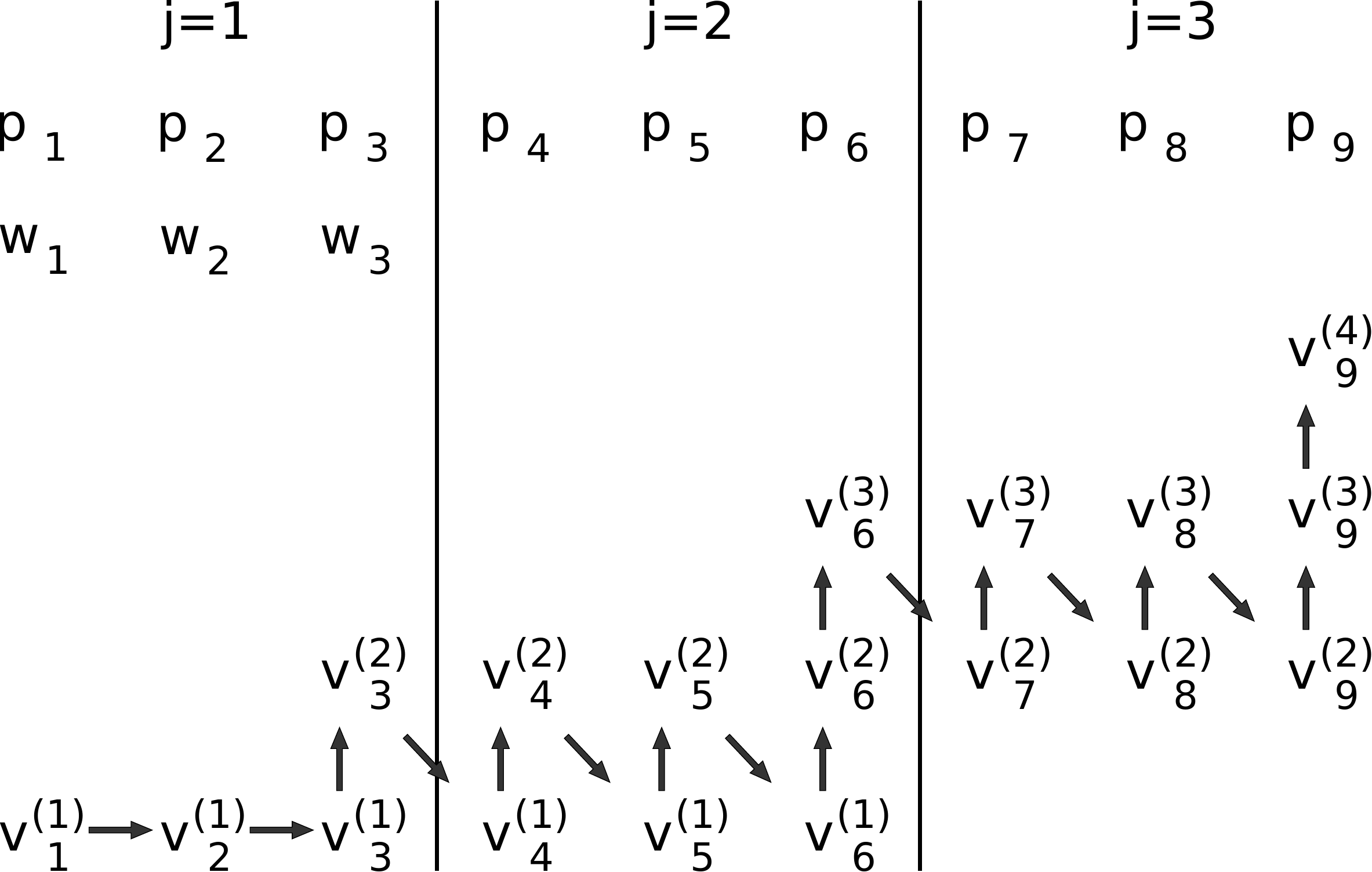}
\caption{
	Order of the computation of the level-vectors for $k=3$. Only vertical and horizontal arrows correspond to a \MV with $\bA$.}
\label{fig:Konstruktion_MLBiCG_arrows}
\end{figure}

\begin{figure}
\centering
\includegraphics[width=1\linewidth]{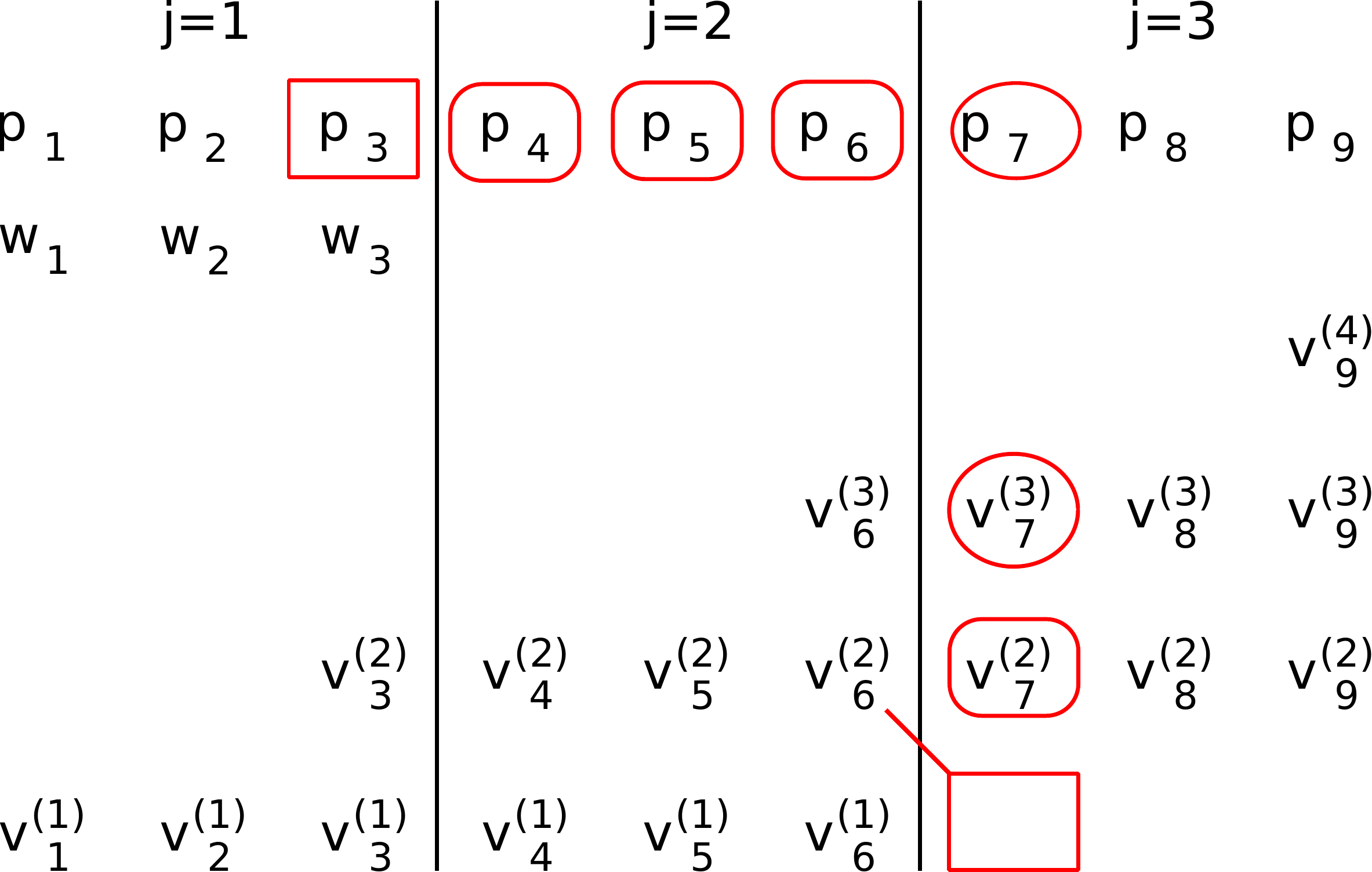}
\caption{
	Example of the orthogonalization procedure for $k=3$ and the $7$-th iteration. Rectangular box: $\bv_6^{(2)} \equiv \bv_7^{(1)}$ is used for computation of the scalar product $\langle \bv_7,\bp_3 \rangle \equiv \langle \bv_7^{(1)},\bw_3 \rangle$. Rounded box:
	Thereafter $\bv_7^{(2)}$ is used to compute the scalarproducts of $\bv_7$ with $\bp_4,...,\bp_6$ by $\langle \bv_7^{(2)},\bw_d \rangle \equiv \langle \bv_7 , \bp_{3+d}\rangle$, $d=1,2,3$. With the resulting linear factors $\bv_7^{(2)}$ is then modified by a scalar combination of $\bv_3^{(2)},...,\bv_6^{(2)}$, so that the vector $\bv_7^{(1)}$ (which is actually not computed) would be orthogonal to $\bp_3,...,\bp_6$. Ellipse: Finally the image of $\bv_7^{(2)}$ with $\bA$ is computed ($\bv_7^{(3)}$) and with this the scalar product $\langle \bv_7,\bp_7 \rangle \equiv \langle \bv_7^{(3)},\bw_1 \rangle$ is evaluated. With the result the vectors $\bv_7^{(2)},\bv_7^{(3)}$ can be scaled, such that $\bp_7^H \cdot \bv_7^{(1)} = 1$ would hold.}
\label{fig:Konstruktion_MLBiCG_recurrsion}
\end{figure}

We write $\bv^{(1)}_i := \bv_i$, $i=1,...,k$. Then we use the following scheme to compute the level-vectors of the following columns of $\bV$:
\begin{algorithmic}[1]
	\For{$j=2,3,...,J$}
	\State $\bv_{j \cdot k}^{(j)} := \bA \cdot \bv_{j \cdot k}^{(j-1)}$
		\For{$i=(j-1) \cdot k,...,j \cdot k - 1$}
			\State \textit{// set of all columns that occur in bi-Lanczos recursion of $(i+1)^{th}$ column}
			\State $\cJ := [\max(1,i-k),...,i] \subset \N$
			\State \textit{// the scalar products for these have to be computed with different level-vectors, therefore define subsets}
			\State $\cJ_{0} := \cJ \cap [(j-1) \cdot k+1,...,j \cdot k] $
			\State $\cJ_{1} := \cJ \cap [(j-2) \cdot k+1,...,(j-1) \cdot k] $
			\State $\cJ_{2} := \cJ \cap [(j-3) \cdot k+1,...,(j-2) \cdot k] $
			\State \textit{// bi-orthogonalizations on level $j-2$}
			\For{all $c \in \cJ_{2}$}
				\State $t_{c,i} := \langle \bw_{r_k(c)} , \bv_i^{(j-1)} \rangle$\quad\textit{// $\equiv \langle \bp_c,\bA \cdot \bv_i \rangle$}
			\EndFor
			\State \textit{// bi-orthogonalizations on level $j-1$}
			\For{all $c \in \cJ_{1}$}
				\State $t_{c,i} := \langle \bw_{r_k(c)} , \bv_i^{(j)} \rangle$\quad\textit{// $\equiv \langle \bp_c,\bA \cdot \bv_i \rangle$}
			\EndFor
			\State $\bv_{i+1}^{(j-1)} := \bv_{i}^{(j)} - \sum_{c \in \cJ_{2} \cup \cJ_{1}} t_{c,i} \cdot \bv_c^{(j-1)}$
			\State $\bv_{i+1}^{(j)} := \bA \cdot \bv_{i+1}^{(j-1)}$
			\State \textit{// bi-orthogonalizations on level $j$}
			\For{all $c \in \cJ_{0}$}
				\State $t_{c,i} := \langle \bw_{r_k(c)} , \bv_{i+1}^{(j)} \rangle$\quad\textit{// $\equiv \langle \bp_c,\bA \cdot \bv_i \rangle$}
			\EndFor
			\State $\bv_{i+1}^{(j-1)} := \bv_{i+1}^{(j-1)} - \sum_{c \in \cJ_{0}} t_{c,i} \cdot \bv_c^{(j-1)}$
			\State $\bv_{i+1}^{(j)} := \bv_{i+1}^{(j)} - \sum_{c \in \cJ_{0}} t_{c,i} \cdot \bv_c^{(j)}$
			\State \textit{// at last we scale}
			\State $t_{i+1,i} := \langle \bw_{r_k(i+1)} , \bv_{i+1}^{(j)} \rangle$\quad\textit{// $\equiv \langle \bp_{i+1},\bv_{i+1} \rangle$}
			\State $\bv_{i+1}^{(j-1)} := 1/t_{i+1,i} \cdot \bv_{i+1}^{(j-1)}$
			\State $\bv_{i+1}^{(j)} := 1/t_{i+1,i} \cdot \bv_{i+1}^{(j)}$
		\EndFor
		\State Delete all $\bv_c^{(d)}$ that were not used in the last walk through the inner for-loop.
	\EndFor
\end{algorithmic}
In the following this scheme is explained. Therefore we also refer to figure \ref{fig:Konstruktion_MLBiCG_arrows}.
\begin{itemize}
	\item Each column of $\bV$ is computed at least in two levels (except the first $k$ columns).
	\item We think of blocks of each $k$ columns. These blocks are numbered by $j=1,...,J$.
	\item For computation of each new column of $\bV$, one needs to compute scalarproducts with shadow vectors $\bp_c$ for indices $c \in \cJ$ from the current and the last two blocks. These scalarproducts can be computed by means of the level-vectors, as with these the scalarproducts with the shadow vectors can be expressed by scalarproducts with the numerical given vectors $\bw_1,...,\bw_k$, cf. \eqref{eqn:LinComb_Lancz}.
	\item By using the level-vectors for the scalarproducts, the linear factors $t_{c,i}$, which are the entries of $\bT$, can be computed.
	\item For computing the linear combination from equation \eqref{eqn:LinComb_Lancz} for the level-vector $(j-1)$ of the current column $i+1$, the level-vectors of the former columns $c \in \cJ_2$ and their linear factors $t_{c,i}$ are used.
	\item For the computation of all the remaining scalarproducts and linear combinations the level of the current column has to be increased during the orthogonalization process. This is the reason why the index set $\cJ$ is subdivided.
	\item The level of a column is increased by just computing its \MV with $\bA$.
	\item At last the lengths of the level-vectors of the new column $i+1$ of $\bV$ are scaled, such that $\bW^H \cdot \bV = \bI$ would hold (although neither $\bW$ nor $\bV$ is actually computed). However these scalings are necessary, so that the entries of $\bT$ fit to the recursion equations, because $\bT$ is used later for the construction of the short representations.
\end{itemize}
The level-vectors $\bv_i^{(j)}$ are computed from left to right, as can be seen in figure \ref{fig:Konstruktion_MLBiCG_arrows}.
The principle of the orthogonalizations by a recursion of length $k+1$ is given in figure \ref{fig:Konstruktion_MLBiCG_recurrsion}:
The frame-types indicate, which level-vector of $\bv_{i+1}$ and which $\bw_\xi$, $\xi=1,...,k$ are combined in a scalar product $\langle \bv_{i+1}^{(j)},\bw_\xi \rangle$ to obtain the value of the scalarproduct $\langle \bv_{i+1},\bp_c\rangle$, $\forall c \in \cJ$. The procedure always takes at first the smallest level-vector (i.e. smallest $j$), computes its scalarproducts, then does the linear combination (which is correlated to the bi-orthogonalization of the original column vector), and then increases the level. Afterwards the scalarproducts with $-$ and orthogonalizations to $-$ the remaining shadow vectors are computed.

\subsubsection{Obtaining the Short Representations} We denote $\bV^\star = [\bv_1^{(1)},...,\bv_k^{(1)},\bv_k^{(2)},\bv_{k+1}^{(2)},...,\bv_{2 \cdot k}^{(2)},\bv_{2 \cdot k}^{(3)},\bv_{2 \cdot k+1}^{(3)},...,\bv_{n-J}^{(J+1)}] \in \C^{N \times n}$ as \textit{level-matrix} of $\bV$.
As a consequence of the above procedure there exists a short recursion for the computation of $\bV^\star$, thus a short representation of $\bV^\star$ can be constructed and stored. The same then holds for its preimage $\bU^\star$, called \textit{level-matrix of $\bU$}. Furthermore, by construction of its columns, $\bV^\star$ can be simply transformed into $\bV$. As consequence one obtains a short representation for $\bV$ over $\bV^\star$.

As the preimage of $\bV^\star$ follows the same recursion as $\bV^\star$, one can also construct a short representation for $\bU$ over a short representation of $\bU^\star$ (where $\bU^\star$ is similarly recaptured such that $\bA \cdot \bU^\star = \bV^\star$).
\largeparbreak
At last we have to clearify how a short representation for $\bW$ can be obtained. For this we store $\tbW = [\bw_1,...,\bw_k]$. We notice that $\hbW = [\tbW,\bA^H \cdot \tbW,...,(\bA^H)^{J-1} \cdot \tbW]$ has the same range as $\bW$, so we now only need a transformation from $\hbW$ to $\bW$.

For this we notice that the linear factors $t_{c,i}$ are just the entries of $\bT$ from the \MLbiL. Thus we can use the recursion
\begin{align*}
	\bA^H \cdot \bW = \obW \cdot \ubT^H\,,
\end{align*}
where $\bW(:,1:k) = \tbW$, and compute the transformation matrix similar to $\bK$ from theorem \ref{thm:ShortReps}.
\largeparbreak
Summing up, with these ingredients the method \SRMLBiCG is now completely described, see remark (calculation procedure). However, we do not postulate a practical implementation of this as we did not already analyse this method that far.

\subsection{A-posteriori-Iterations for \SRMLBiCG}

\paragraph{Motivation}
After \SRMLBiCG has been used to compute a recycling solution $\bx$, its residual $\br$ is orthogonal to a special structured test space: $\br \perp \cP = \cK^\star_J(\bA;[\bw_1,...,\bw_k])$.

As we already did in chapter \ref{chap:ApostImprov}, we want to derive a scheme for a-posteriori-iterations for a residual and recycling data obtained from \SRMLBiCG, in which a property comparable to $\br \perp \cP$ does not get lost.

In chapter \ref{chap:Trans_K_to_G} we achieved this for the case $s \geq k$ where $k = 1$. Now we want to consider the case where $k > 1$ because only for this choice \SRMLBiCG is more efficient than \SRBiCG in terms of \rMV for the solution of a first system.

\paragraph{\IDR for \SRMLBiCG}
After computation of the \MLbiL and computation of a recycling solution we are given data $\br,\bv_{n+1} \perp \cP$. (Actually we do not have $\bv_{n+1}$ explicitly but let us assume for this derivation that we would have it.) Then by just expanding the size of the \MLbiL from size $n$ to size $n+s$, $s=k$ we can get $k$ columns $\bv_{n+1},...,\bv_{n+k} \perp \cP$, instead of only $\bv_{n+1}$, where $\cP = \opImage(\bW)$ and $\bW$ is the shadow basis from the \MLbiL of size $n$, thus $\bW \in \C^{N \times n}$.

In principle we could go ahead as in chapter \ref{chap:Trans_K_to_G}. Instead of only using $\bv_{n+1}$ as in chapter \ref{chap:Trans_K_to_G}, we would now use $\bv_{n+1},...,\bv_{n+k}$. And instead of \textit{one} $\hbw = \bw_{n+1}$ as in chapter \ref{chap:Trans_K_to_G} we would now need $\hbw_1 = \bp_{n-k+1}$,$\hbw_2 = \bp_{n-k+2}$,...,$\hbw_k = \bp_{n}$ to ensure that the residual will remain orthogonal to $\cP$. Using these data one obtains exactly a generalization of the approach from chapter \ref{chap:Trans_K_to_G}.

However, this approach has the problem that we barely got rid of computing the columns $\bp_1,...,\bp_n$, so we do not want to compute them only to enable a-posteriori-orthogonality.

That is why in the case of \SRMLBiCG I propose to proceed \textit{not} as in chapter \ref{chap:Trans_K_to_G}, but instead using equations like \eqref{eqn:IDR_t} and \eqref{eqn:IDR_g}. By this equations I theoretically derive a method as follows:

Eqs. \eqref{eqn:IDR_t}, \eqref{eqn:IDR_g} mean multiplying both residual $\br$ and auxiliary vectors $\bv_{n+1},...,\bv_{n+k}$ by a polynomial $\Omega_J(\bA)$ of degree $J$ with non-zero leading coefficient, such that one obtains vectors $\bt$ and $\bg_1,...,\bg_k$  (notation analogue to chap. \ref{chap:Trans_K_to_G}). $\bt$ then is the new residual for which $\bx$ has to be recaptured, and $\bg_i$,$i=1,..,k$ are the new auxiliary vectors for which the preimages also have to be hold.

During the computation of the \MLbiL we compute level-vectors of the columns $\bv_i$ of the matrix $\bV$ up to level $J+1$. Thus for the choice of $\Omega_J(\bA) = (\bI-\rho \cdot \bA)^J$ we already have computed $\bg_i$, $i=1,..,k$ within our procedure for computation of the \MLbiL due to the computation of the $\bv_i$'s level-vectors: We have $\bg_i = \bv^{(J+1)}_{n+i} \equiv \bA^J \cdot \bv_{n+i}$,$i=1,..,k$ (we can ignore the shift as a shift does not change the basis vectors of a Hessenberg-decomposition). So we get the auxiliary vectors for free. Unfortunately for the residual we then have to use its update $\bt := (\bI-\rho \cdot \bA)^J \cdot \br$, which is only good if there exists an appropriate $\rho$ or if $\bI-\rho \cdot \bA$ can be replaced by an iteration matrix of a convergent basic iterative scheme for the according problem. Otherwise and for large choice of $J$ this is cursed to fail.

To clearify this: At this point and at this moment I do not have a better idea. (In fact one could also choose $\Omega$ of degree lower then $J$ and by this just looses a few levels of the Sonneveld-space in which the residual is actually moved.)

\section{Preconditioning}

In this chapter it is shown how one can apply preconditioning for \SRIDR, \SRBiCG, \SRCG and \SRMR.

As a restriction I do not focus on flexible, changing or non-linear preconditioners. These restrictions are reasonable for short-recurrence methods. After these restrictions it is well-known how to apply preconditioning strategies.
\largeparbreak
From a linear system
\begin{align*}
\bA \cdot \bx = \bb
\end{align*}
by substitutions we end up with
\begin{align*}
	\bL^{-1} \cdot \bA \cdot \bR^{-1} \cdot \tbx &= \tbb\,,
\end{align*}
where we aim for solution after $\tbx$. In this $\bM = \bL \cdot \bR \approx \bA$, where $\bL,\bR \in \C^{N \times N}$ are just two matrices with full rank.

\paragraph{Splitted Preconditioning} In case of $\bL,\bR \in \C^{N \times N}$ given, applications with $\bA$ in \SRIDR, \SRBiCG (and, assumed $\bR = \bL^H$ and $\bA = \bA^H$, also in \SRMR) can be simply replaced by products with $\bL^{-1} \cdot \bA \cdot \bR^{-1}$.

\paragraph{Unsplitted Preconditioning} Another situation is, where only applications of matrix-vector-products with $\bM^{-1}$ and $\bA$ can be performed. In the unsymmetric case, one could solve
\begin{align*}
\bM^{-1} \cdot \bA \cdot \bx = \bM^{-1} \cdot \bb\,.
\end{align*}
For this end, one would replace matrix-vector-products with $\bA$ in \SRIDR and \SRBiCG by matrix-vector-products with $\bM^{-1} \cdot \bA$.

\paragraph{Unsplitted Preconditioning $-$ Preserving Symmetry} One special case in this is, if $\bA$ is Hermitian and $\bR = \bL^H$, i.e. $\bM$ is positive definite. In this case one can still apply \SRMR and \SRCG, although matrix-vector-products with the Hermitian matrix $\bL^{-1} \cdot \bA \cdot \bL^{-H}$ cannot be computed, but instead only products with $\bA$ and $\bM^{-1}$.

Methods based on the symmetric Lanczos-procedure, that are capable of such kind of preconditioning, are well-known, e.g. preconditioned \CG and \MINRES, cf. \cite[chap. 9, algo. 9.1]{Saad1}, \cite{Sym-Pres-MINRES}. There also exist symmetry-preserving \GMRES-variants \cite{Sym-Pres-GMRES}. As in this work we need matrices for which we can store short representations, I briefly review the conceptual idea of these.

Instead of computing a Lanczos-decomposition of the form
\begin{align*}
	\bL^{-1} \cdot \bA \cdot \bL^{-H} \cdot \bV = \obV \cdot \obT\,,
\end{align*}
with $\obT = \operatorname{tridiag}(\beta_2,...,\beta_{n+1};\alpha_1,...,\alpha_n;\beta_2,...,\beta_n) \in \C^{(n+1) \times n}$, one uses a transformed basis matrix $\obZ = \bL^{-H} \cdot \obV \in \C^{N \times (n+1)}$. From this it follows
\begin{align*}
	\bM^{-1} \cdot \bA \cdot \bZ = \obZ\cdot \obT\,.
\end{align*}
The tridiagonal matrix $\obT$ is not affected by this basis transformation as
\begin{align*}
	\obT = \obV^H \cdot \bL^{-1} \cdot \bA \cdot \bL^{-H} \cdot \bV = \obZ^H \cdot \bA \cdot \bZ\,.
\end{align*}
From the latter equation we see that in the $i^\text{th}$ iteration, i.e. the computation of the $(i+1)^\text{th}$ column of $\bZ$ from the former columns, the entries $\alpha_i$ and $\beta_{i+1}$ of $\obT$ can be already computed, thus the linear-factors for the next recursion
\begin{align*}
\bz_{i+1} = 1/\beta_{i+1} \cdot (\bM^{-1} \cdot \bA \cdot \bz_i - \alpha_i \cdot \bz_i - \beta_i \cdot \bz_{i-1})
\end{align*}
are known. So one can now store a short representation of $\bZ$ (in the V-approach) or of its preimage $\bF$ (with $\bZ = \bM^{-1} \cdot \bA \cdot \bF$, in the U-approach).

A formula for computation of recycling solutions by means of the short representation can be simply derived, e.g. for the V-approach we have
\begin{align*}
	\tbx &= \bV \cdot \obT^\dagger \cdot \obV^H \cdot \tbb\\
	\bx  &= \bL^{-H} \cdot \bV \cdot \obT^\dagger \cdot \obV^H \cdot \bL^{-1} \cdot \bb\\
	     &= \bZ \cdot \obT^\dagger \cdot \obZ^H \cdot \bb\,.
\end{align*}

\paragraph{Remarks}
\begin{itemize}
	\item We have implementations of \SRIDR and \SRBiCG that enable preconditioning, as the routine for matrix-vector-products with $\bA$ can be replaced by a user-defined function for e.g. products with $\bL^{-1} \cdot \bA \cdot \bR^{-1}$.
	\item For \SRMR we have implemented the unsplitted symmetry-preserving preconditioning approach.
\end{itemize}

\section{Numerical Experiments}\label{chap:Experiments}
The following experiments shall indicate the potential of the proposed principles and thereon based methods, give evidence for the theoretical findings discussed above, and wake the reader's interest on these methods.
\largeparbreak
For the experiments we constructed very plain testcodes for \RGCR, \SRIDR and \SRBiCG.

As the numerical residuals decouple from the approximate solutions in \SRIDR due to peaks in the convergence graphs, we used cyclic residual replacements.

For \SRBiCG we implemented the stabilization-by-blocking approach and the a-posteriori-iterations with orthogonality-conservation. For simplicity we used a uniform block size for the stabilization and \IDReins for the a-posteriori-iterations.

In the figures below we presented the \enquote{true} residuals, i.e. we plotted $\|\bb-\bA \cdot \bx\|$ for the numerical values of $\bx$.

Further information on the testcodes and accessibility can be found on \texttt{MartinNeuenhofen.de} $\rightarrow$ \texttt{Short Recycling} .

\paragraph{Experimental setting}
There are three test cases, each of which is solved without preconditioning. For each case the \rhs-es are solved \underline{one after another}.
\begin{Aufz}
	\item The ocean-problem form \cite{IDRweb} with \texttt{stommel\_4.mtx} and twelve \rhs-es; $N = 2594$, $\cond_1(\bA) \approx 5.1858 \cdot 10^5$.
	\item A convection-diffusion-reaction-problem on the unit cube from \cite{IDRweb} with standard parameters except $h=0.05$ and the provided \rhs from the test routine; $N = 6858$, $\cond_1(\bA) \approx 336.48$.
	\item Matlab's Poisson-matrix \cite{MATLAB-TEST} with $100 \times 100$ inner points and \rhs \texttt{ones(10000,1)}; $N = 10000$, $\cond_1(\bA) \approx 6.0107\cdot 10^3$.
\end{Aufz}
For the second and third test case so far we only have one \rhs. From this we construct $z=5$ \rhs-es for the second test case and $z=10$ \rhs-es for the third test case by constructing a matrix $\bB \in \C^{N \times z}$ with \rhs-es as columns. These are constructed as follows:
\begin{algorithmic}[1]
\State $\bB(:,1) := \bb$
\For{$i=1,...,z$}
	\State $\bB(:,i+1) := \bA \backslash \bB(:,i)$
	\State $[\bB,\sim] := $\texttt{qr($\bB,0$)} \quad \textit{// orthogonal matrix of reduced QR-decomposition}
\EndFor
\end{algorithmic}
This is due to the following reasons:
\begin{Aufz}
	\item The \rhs-es are orthogonalized w.r.t. each orther as otherwise the former solutions could be used as initial guess to ensure orthogonality of the remaining residual.
	\item The space of \rhs-es is build as a reverse Krylov-subspace as such structures e.g. occur in the introductory example \eqref{eqn:Fourier} just like in many other discretizations with implicit time integration schemes. In case of finite element discretizations there is still a mass-matrix, but for simplicity and as this is often well-conditioned we neglect it here.
\end{Aufz}

\paragraph{Proposition (foreward and reverse Krylov-subspaces)}
We notice that for well-conditioned positive definite systems the Krylov-subspace is always a very usefull recycling space for the solution to a \rhs in a backward Krylov-subspace $\cK_p(\bA^{-1};\bb\hia)$ of small degree $p$.

\subsection{Ocean-Problem}

\paragraph{\RGCR}
At first we validate that recycling is usefull at all. Therefore we solve the problem first with full \RGCR. For this we store all basis vectors that were computed in the first solution process and reuse this data for solution of all the following systems.

\begin{figure}
\centering
\includegraphics[width=0.7\linewidth]{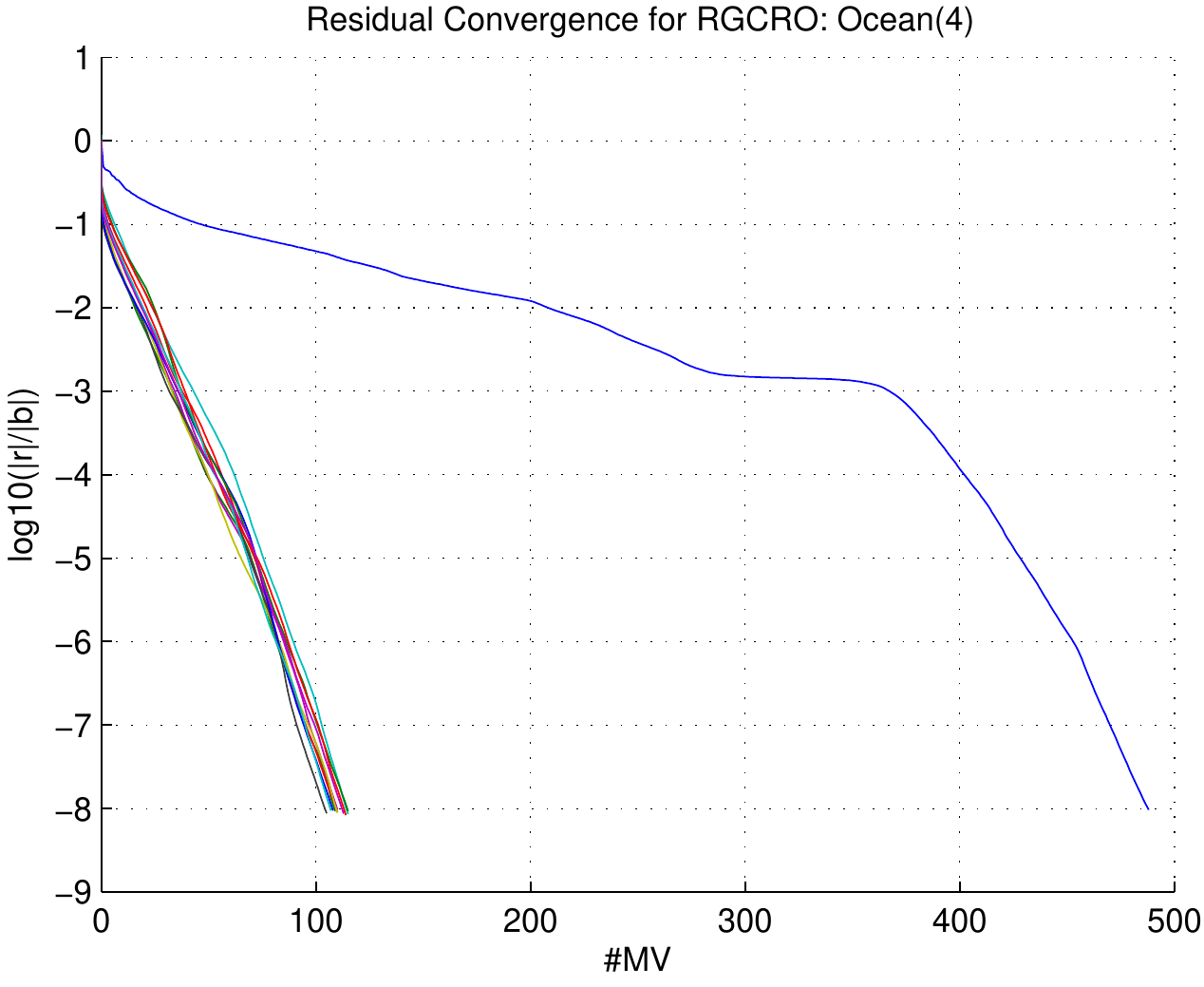}
\caption{\RGCR for Ocean.}
\label{fig:Ocean4_GCR}
\end{figure}

In figure \ref{fig:Ocean4_GCR} one can see that \RGCR needs roughly $500$ \MV-s for solution of the first system. By means of recycling then the solutions to the following $11$ \rhs-es can be found each in rougly $100$ \MV-s. So for this problem recycling-methods in fact turn out to be very usefull!

\begin{figure}
	\centering
	\includegraphics[width=0.7\linewidth]{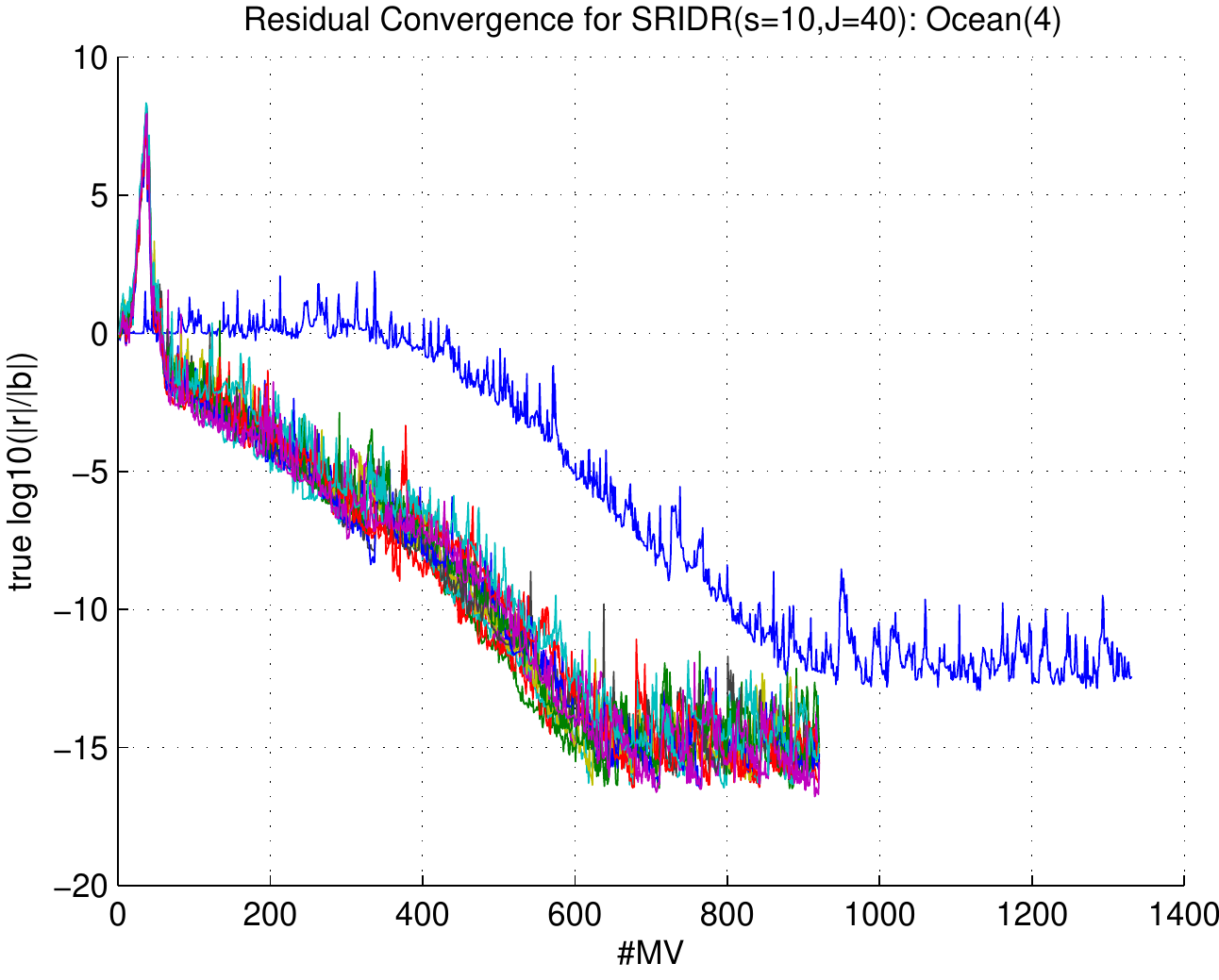}
	\caption{\SRIDR for Ocean.}
	\label{fig:Ocean4_SRIDR}
\end{figure}

\paragraph{\SRIDR}
In figure \ref{fig:Ocean4_SRIDR} the same problem is solved with \SRIDR for $s=10$. The data for the recycling was fetched after $J^\star = 40$ \IDRO-cycles of the first system, thus the recycled test space consists of $400$ dimensions. I observe and interprete the convergence graphs as follows:
\begin{Aufz}
	\item For the first system my \IDR implementation performs a bit worse than the \IDR-biortho Matlab implementation of van Gijzen, cf. \cite{IDRweb}. This can have several reasons, e.g. I added a cyclic residual replacement.

	\item The peaks at the beginning of the recycling processes probably come from the lack of a residual minimizing property during the recycling-\IDRO-cycles.

	\item Up to relative residual tolerance $\tol = 10^{-10}$ one can compute recycling solutions in $2/3$ of the \rMV-s compared to the solution process of the first system that is without recycling.

	\item The relative residual can be reduced to $10^{-15}$ by use of residual replacement.

	\item However the performance of \SRIDR does still not reach the performance of full \RGCR.
\end{Aufz}

\paragraph{\SRBiCG}
\SRBiCG fails for this test case. The reason for this is, that my \BiCG implementation does not converge for this problem (in contrast to Matlab's \texttt{bicg}, which still converges in $1200$ \MV-s).

\subsection{Convection-Diffusion-Reaction}
\begin{figure}
\centering
\includegraphics[width=0.8\linewidth]{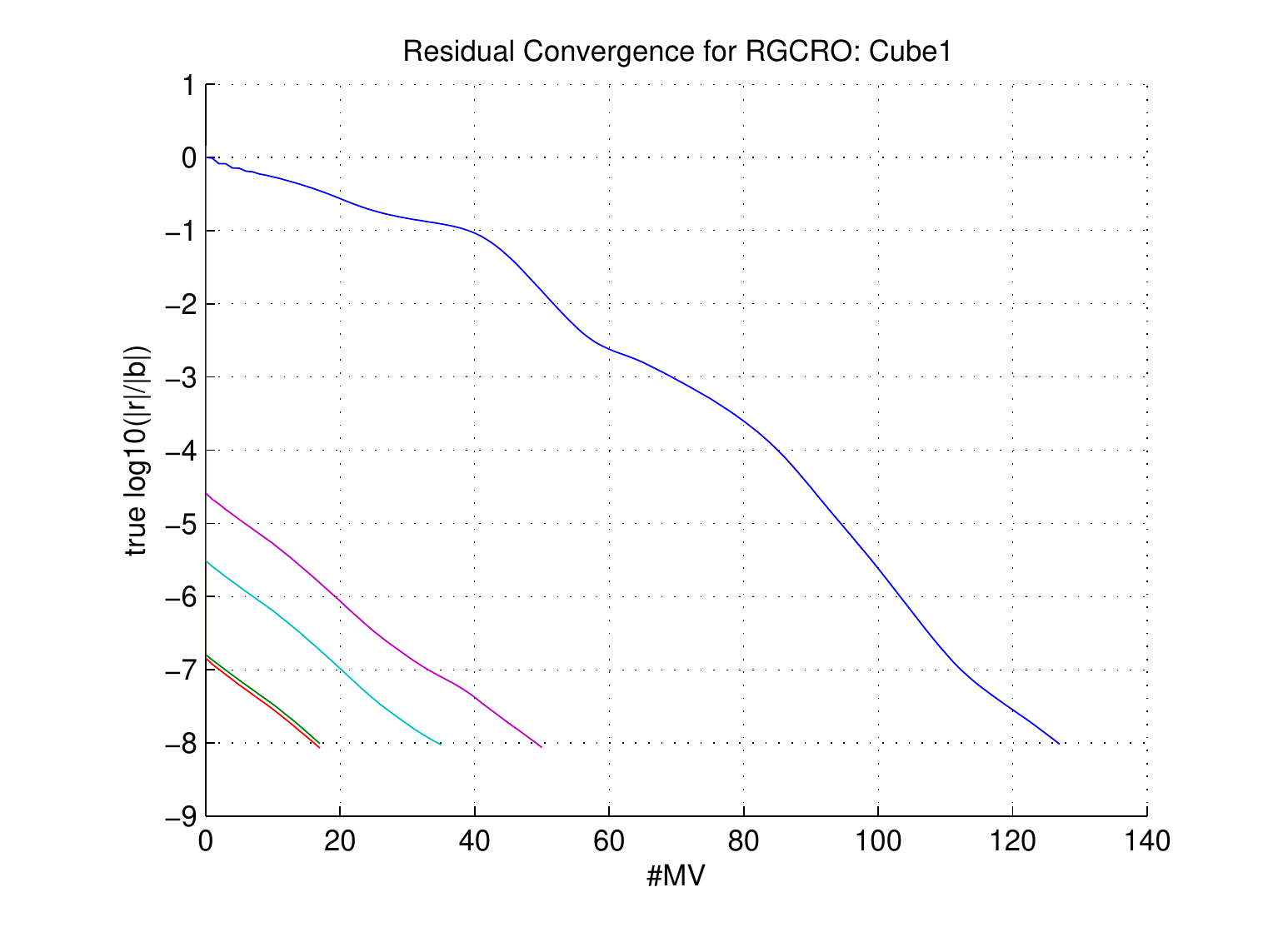}
\caption{\RGCR for Convection-Diffusion-Reaction.}
\label{fig:Cube1_GCR}
\end{figure}

\paragraph{\RGCR}
As before, we validate with \RGCR, that recycling brings some benefit for this problem instance. As before, we store all basis vectors from the solution process of the first \rhs and reuse them for the following \rhs-es. The convergence curves are given in figure \ref{fig:Cube1_GCR}. One can see, that for the second and third \rhs the solution lies nearly in the recycled Krylov-subspace, but for the fourth and fifth \rhs \RGCR needs significantly more additional \MV-s. Also for this test case recycling enables a large cost reduction in \rMV.

\begin{figure}
	\centering
	\includegraphics[width=0.8\linewidth]{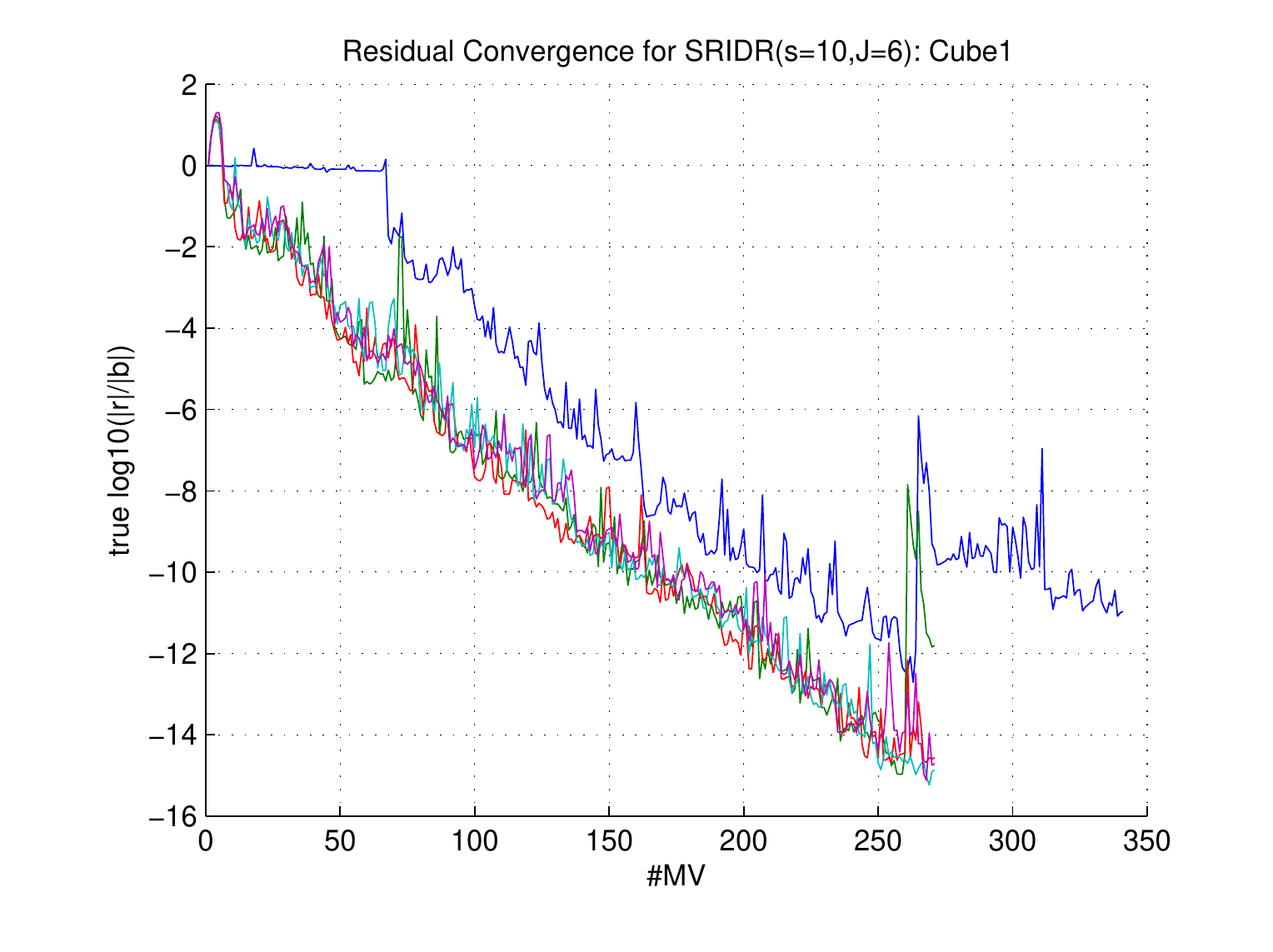}
	\caption{\SRIDR for Convection-Diffusion-Reaction.}
	\label{fig:Cube1_SRIDR}
\end{figure}

\paragraph{\SRIDR}
For \SRIDR we used $s=10$ and fetched recycling data after $J^\star = 6$ \IDRO-cycles in the solution process for $\bb_1$. Thus the recycled test space has only $60$ dimensions, which is only half the size of the full \RGCR space. I chose $J^\star$ that small for the following reasons:
\begin{Aufz}
	\item Due to \SRIDR's use of another search space for the recycling one cannot expect that this search space offers the possibility to provide a recycling solution which is as good as the one of full \RGCR, cf. spaces of \SRIDR, chap. \ref{chap:SRIDR-Spaces}.

	\item Instead \textit{I assume} that the best behaviour that one can expect is that the convergence curve starts right at the beginning with the slope from the point of that curve from where the recycling data was obtained.
\end{Aufz}
In figure \ref{fig:Cube1_SRIDR} the convergence curves are given in blue for the first system and in other colors for the subsequent \rhs-es. One can see that by use of the recycling the primary stagnation phase can be skipped.

For a relative residual tolerance of $10^{-10}$ the recycling solutions can be obtained with approximately $50$ fewer \MV-s then the first solution. But as before, the cost reduction in \rMV for \SRIDR is far below the one of full \RGCR.

\begin{figure}
	\centering
	\includegraphics[width=0.8\linewidth]{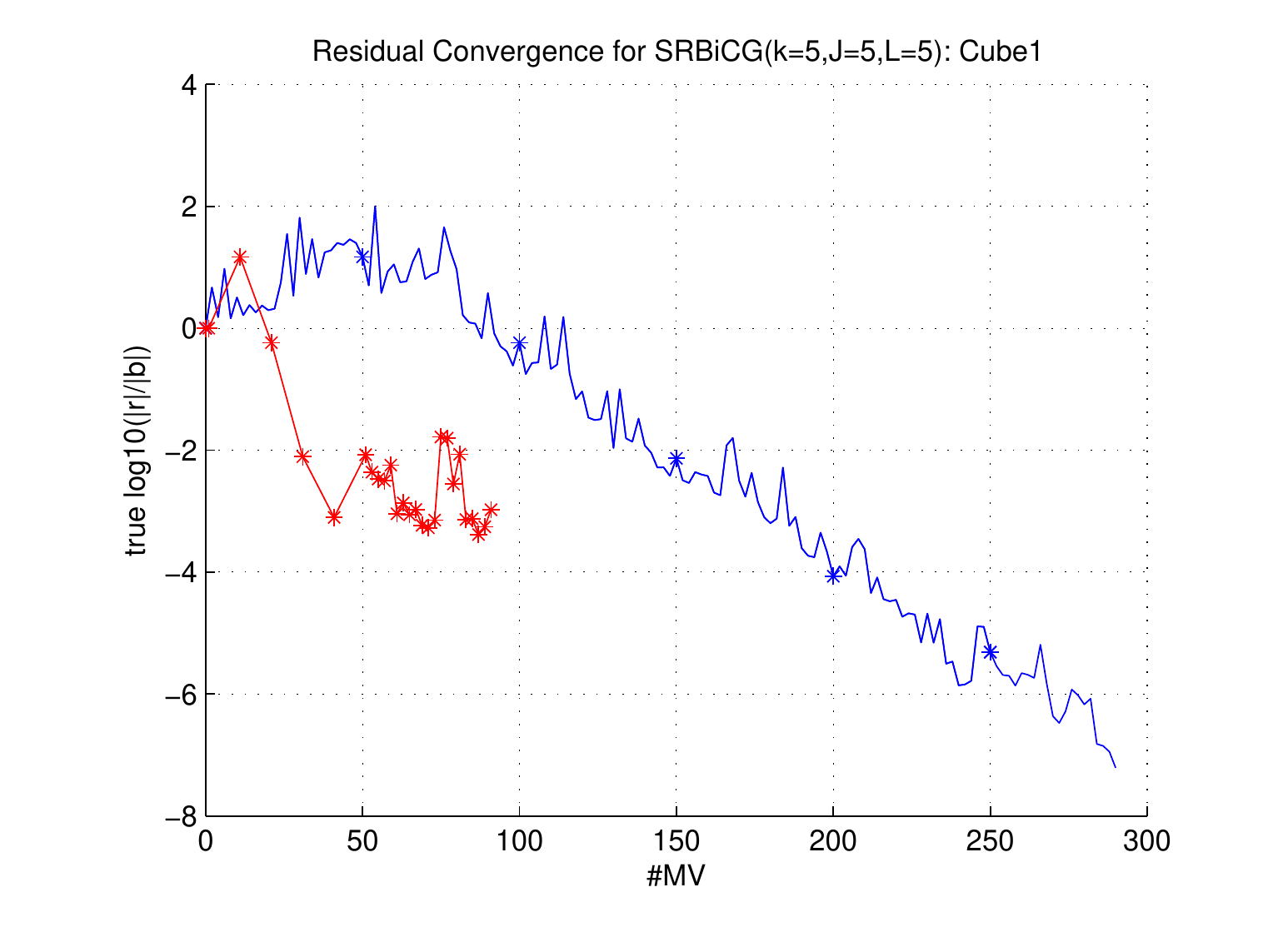}
	\caption{\SRBiCG for Convection-Diffusion-Reaction.}
	\label{fig:Cube1_SRBiCG}
\end{figure}

\paragraph{\SRBiCG}
\SRBiCG also fails for this test problem. To find a possible reason to this, we solved two consecutive times for the first \rhs: Once without recycling (blue), and afterwards once with recycling (red) by use of the data obtained from the first solution process, cf. fig. \ref{fig:Cube1_SRBiCG}.

The chosen parameters are $\ell = 5$ blocks of each $k=5$ stored columns (thus overall for both search- and test-space each $\ell \cdot k = 25$ columns have to be stored) and $J=5$. By this the respective blocks have  $k \cdot J = 25$ columns. All in all a $\ell \cdot k \cdot J = 125$-dimensional search- und test-space are recycled.

We marked the respective blocks in the convergence graphs by \texttt{*}-marks. Under exact arithmetic the marks of the red and blue graph would be pairwise on the same height, as the \rhs is orthogonalized in both cases with the same search space on the same test space (in blue by the original \BiCG process and in red by the short recycling procedure).
For the first $3$ blocks (the relative residual then is at ca. $10^{-2}$) the short recycling procedure seems to work. However for the fouth block the scheme fails, as can be seen from the large gap between the red and the blue mark.

A possible reason for this is that either the data-accuracy or the representation-accuracy of the fourth block is bad and that this leads to the failure in this short recycling step.

A remedy to this could be to choose the block size adaptively, e.g. by the condition number of the current tridiagonal block $\bT^{(i)}$ (cf. notation from chap. \ref{chap:Robust}).

The following a-posteriori-iterations (beginning in the diagram in fig. \ref{fig:Cube1_SRBiCG} at \rMV$>50$) fail then because they use descent directions orthogonal on the columns of $\bW$ which already have not been eliminated properly from the residual.

\subsection{Poisson}
\begin{figure}
\centering
\includegraphics[width=0.8\linewidth]{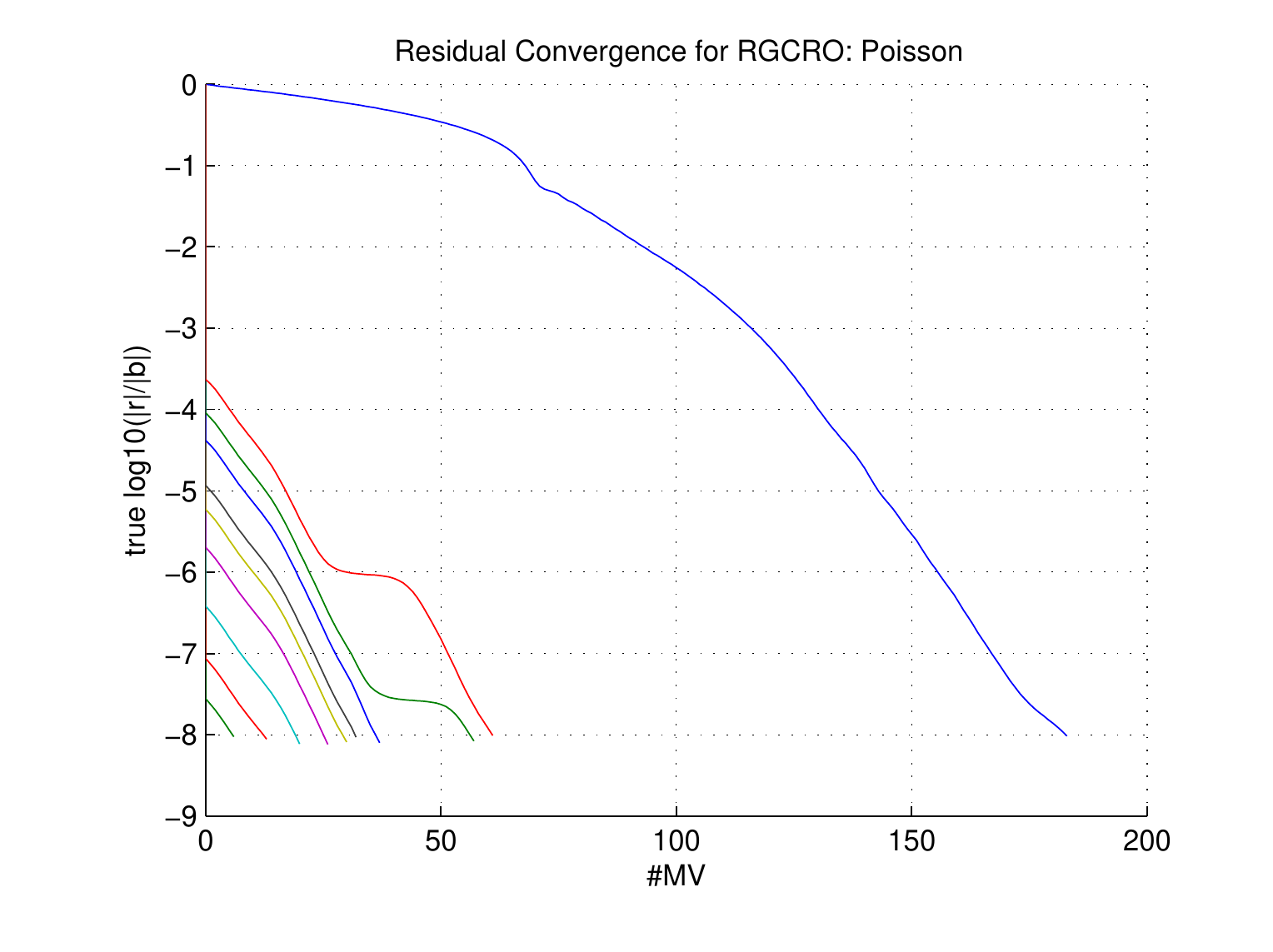}
\caption{\RGCR for Poisson.}
\label{fig:poisson_GCR}
\end{figure}

\paragraph{\RGCR}
As before we check the usefullness of recycling procedures for this problem by recycling all data that is obtained from the solution process for the first \rhs. The convergence graphs are given in figure \ref{fig:poisson_GCR}.

Similar to the convection-diffusion-reaction test case, the use of recycling decreases for subsequent \rhs-es but remains worthwhile up tp the last system.

\paragraph{\SRIDR}
For this we fetch recycling data after the $J^\star=15$-th \IDRO-cycle for $s=10$ of the first system. This is approximately after $150$ \MV-s. The recycled test space has dimension $s \cdot J^\star = 150$. This is little below the dimension of the the recycled test- and search-space of \RGCR from the first system (about $180$ dimensions for \RGCR).

Comparable to the result for the convection-diffusion-reaction problem one can approximately skip the first $80$ \MV-s from the first (blue) convergence graph (where the convergence stagnates) for the solutions obtained by recycling. A comparable performance to that can also be obtained by only choosing $J^\star=8$.

\begin{figure}
	\centering
	\includegraphics[width=0.8\linewidth]{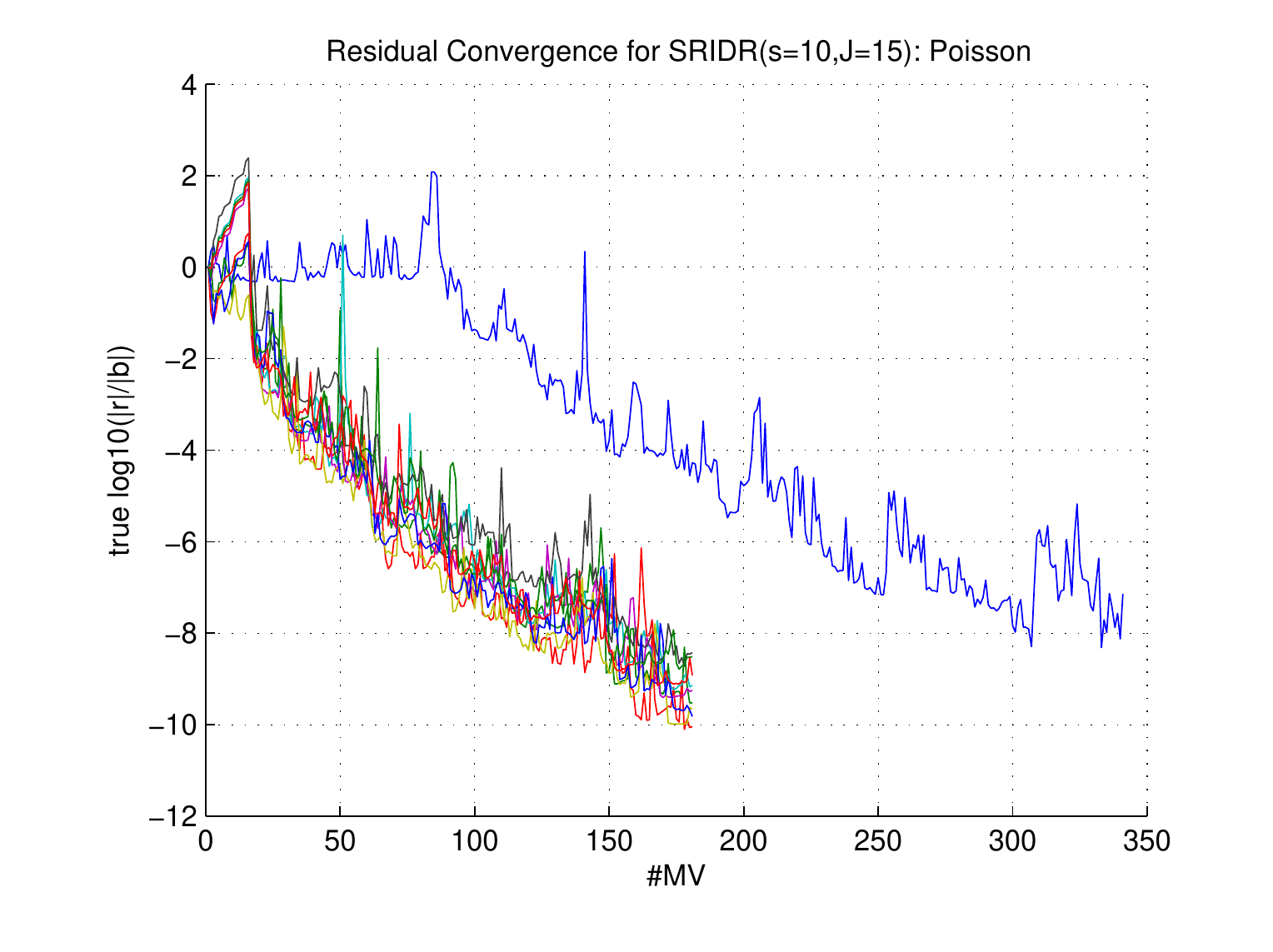}
	\caption{\SRIDR for Poisson.}
	\label{fig:poisson_SRIDR}
\end{figure}

\paragraph{\SRBiCG}
As this problem is symmetric and positive definite, \SRBiCG just performes like \CG (because $\bw_1 = \bu_1 \parallel \bb_1$ is chosen). Therefore good convergence at least for the first system should be guaranteed.

\begin{figure}
\centering
\includegraphics[width=0.8\linewidth]{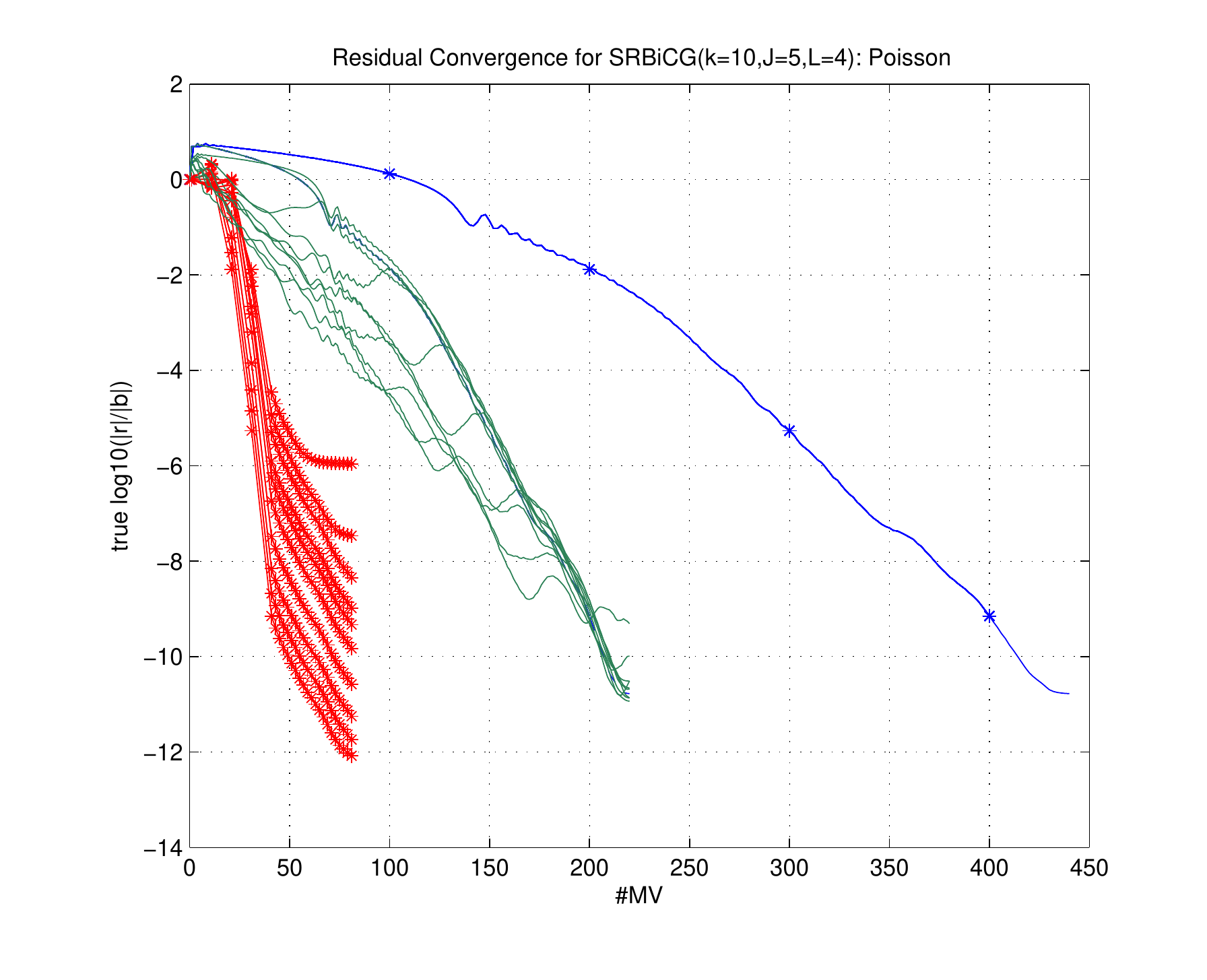}
\caption{\SRBiCG for Poisson.}
\label{fig:poisson_SRBiCG}
\end{figure}

In figure \ref{fig:poisson_SRBiCG} one can see the following convergence graphs: In blue the convergence graph of \BiCG for the first system is shown. The recycling solutions for the first to the tenth system were computed with recycling data from the solution process which is represented in the blue curve. The corresponding curves are depicted in red. Independent of the \SRBiCG-computations in green the convergence graph of \CG is plotted for each \rhs.

We see that \BiCG needs double as many \rMV-s as \CG. As already noticed instead of using \SRBiCG for this test case, one could also use \SRCG, which would produce the same red curves and a blue curve with halved scaled horizontal axis.

For the recycling I used a $200 = \ell \cdot k \cdot J$-dimensional search- and test-space, divided into $\ell = 4$ blocks of each $k=10$ stored columns $\tbU^{(i)},\tbW^{(i)} \in \C^{N \times k}$, $i=1,...,\ell$ for both search- and test-space, and $J=5$ for each block.
Thus the overall number of stored columns is $\ell \cdot k = 40$ columns for each search- and test-space (notice that as $\tbU^{(i)}=\tbW^{(i)}$, one could halve the number of stored columns in this case due to $\bA^H = \bA$).

One sees, that by use of recycling the residuals to the solutions of the \rhs-es $\bb_1,...,\bb_4$ can be reduced below relative size of $10^{-8}$ within $50$ \MV-s. This is still comparable to full \RGCR!

\section{Conclusion}
Essentially I have introduced two principles in this work, that can be used to recycle information from an iterative solution process of a linear equation system, to enable faster solution of subsequent \rhs-es to the same system matrix. These principles are based on short recurrences and small memory requirements.
\largeparbreak
The first of these principles works on Sonneveld-spaces and offers an elegant theory providing large theoretical insight. It is simple to implement and enables many options for modifications and extensions such as the a-posteriori-orthogonality-conservation. However in practice the naive \SRIDR-method without sophisticated choice of relaxations offers rather moderate speed-ups, comparable to methods based on compressing, by just keeping the best convergence slope from the former solution process. Compared to a full recycling method like \RGCR, the naive \SRIDR implementation that was tested here performs not so well.

\largeparbreak
The other principle is less elegant but offers the ability for recycling of the original search space. It copies the algorithmic trick of \SRIDR of using a Horner scheme for construction of a search space candidate by only storing a few columns. This principle leads to block-Krylov-matrices, embedded by the theorem of short representations. It is more complicated to implement, is error-prone due to round-off and complicated basis construction schemes, and it does not offer any new theoretical insights. It is just a geometrical construction approach.

In fact every theoretical supplement to the methods of this principle throughout this text has been derived from the \IDRO-principle (e.g. Horner scheme, a-posteriori-orthogonality).

As an advantage on the other hand - if it does not fail - this approach of short representations seems to be more efficient then \SRIDR as the original Krylov-subspace is used as search-space for the recycling solutions, cf. figure \ref{fig:poisson_SRIDR} compared to \ref{fig:poisson_SRBiCG}.

A second advantage of this more geometrical approach is that at least we know what has to be stabilized. Thus for the short recycling approaches we can develop stabilization techniques, whereas for \SRIDR we even do neither know what may be need to be stabilized nor what numerical property would be desirable, nor what stability measure could be even applied.

\largeparbreak
These topics may be considered for future work:
\begin{Aufz}
	\item reliable implementations,
	\item reliable adaptive stabilization-by-blocking approaches,
	\item finding other translations from \SRIDR to short recycling for conventional, Krylov-subspace methods,
	\item study of reverse Krylov-subspaces; for which cases is a recycling method probably advantageous?,
	\item solving sequences with changing matrices,
	\begin{Aufz}
		\item using \SRIDR as interior linear method of an outer Krylov-subspace method,
		\item using \SRIDR as it is for changing matrices,
	\end{Aufz}
	\item numerical experiments with preconditioning.
\end{Aufz}

\section{Acknowledgements}
I would like to thank Peter Sonneveld and Martin van Gijzen for their positive feedback, interest and helpfull discussions on several methods described in this text. For me this was the driving factor for this work.

I also would like to thank Arnold Reusken for his support of my scientific work and his help to get me into contact with experts on the methods I deal with.

My special thanks go to Sven Groß, for his time investment in corrections of my drafts, for long and helpfull discussions, learnings on scientific writing, and for his personal motivation.

\largeparbreak
\largeparbreak
\largeparbreak
\FloatBarrier

\section{References}
\begingroup
\renewcommand{\section}[2]{}%

\endgroup
\end{document}